\documentclass[11pt]{amsart}
\usepackage[utf8]{inputenc}
\usepackage{mathpazo}
\usepackage{hyperref}
\usepackage[margin=1in]{geometry}
\usepackage{array,multirow}
\usepackage[alphabetic,backrefs,lite]{amsrefs}
\usepackage[all]{xy}
\usepackage{array}
\usepackage{amssymb,amsmath,latexsym,amsthm}
\usepackage{tikz-cd}
\usepackage{cleveref}
\usepackage{comment}
\usepackage{mathtools}
\usepackage{extpfeil}

\newtheorem{theorem}{Theorem}[section]
\newtheorem{lemma}[theorem]{Lemma}
\newtheorem{proposition}[theorem]{Proposition}
\newtheorem*{theorem*}{Theorem}

\newtheorem{assumption}[theorem]{Assumption}

\newtheorem{corollary}[theorem]{Corollary}
\newtheorem{predefinition}[theorem]{Definition}
\newenvironment{definition}{\begin{predefinition}\rm}{\end{predefinition}}
\newtheorem{preremark}[theorem]{Remark}
\newenvironment{remark}{\begin{preremark}\rm}{\end{preremark}}
\newtheorem{prenotation}[theorem]{Notation}
\newenvironment{notation}{\begin{prenotation}\rm}{\end{prenotation}}
\newtheorem{preexample}[theorem]{Example}
\newenvironment{example}{\begin{preexample}\rm}{\end{preexample}}

\numberwithin{equation}{section}

\newcommand \conn {{\rm conn}}
\newcommand \ZZ {{\mathbb Z}}
\newcommand \QQ {{\mathbb Q}}
\newcommand{\bQ}{{\mathbb Q}}
\newcommand \CC {{\mathbb C}}
\newcommand \RR {{\mathbb R}}
\newcommand \FF {{\mathbb F}}

\newcommand{\Gl}{G_{A,\ell}}
\newcommand{\Glcon}{G_{A,\ell}^{\circ}}

\newcommand{\adeles}{{\mathbb A}}
\newcommand{\bG}{{\mathbb G}}

\newcommand \cf {{\mathfrak f}}
\newcommand \sss {{\rm ss}}
\newcommand \ord {{{\rm ord}}}
\newcommand \tr {{\rm tr}}
\newcommand \Fr {{\rm Frob}}

\DeclareMathOperator{\End}{End}
\DeclareMathOperator{\Hom}{Hom}
\DeclareMathOperator{\Gal}{Gal}
\DeclareMathOperator{\GSp}{GSp}
\DeclareMathOperator{\Sp}{Sp}
\DeclareMathOperator{\GL}{GL}
\DeclareMathOperator{\SL}{SL}

\DeclareMathOperator{\Lie}{Lie}
\DeclareMathOperator{\Res}{Res}
\DeclareMathOperator{\der}{der}
\DeclareMathOperator{\id}{id}
\DeclareMathOperator{\diag}{diag}
\DeclareMathOperator{\spn}{span}

\newcommand \cO {{\mathcal O}}
\newcommand \cU {{\mathcal U}}
\newcommand \vvv {{\mathfrak{p}}}
\newcommand{\fp}{{\mathfrak{p}}}
\newcommand{\fP}{{\mathfrak{P}}}
\newcommand{\fS}{{\mathfrak{S}}}
 
\newcommand \cB {{\mathcal B}} 
\newcommand \cP {{\mathcal P}} 

\newcommand{\bbb}{T}
\newcommand \att {t}
\newcommand \rr {\kappa}
\newcommand{\Wc} {{W^\circ}}
\newcommand{\Lbar}{{\overline{L}}}
\newcommand{\Qbar}{{\overline{\mathbb{Q}}}}
\newcommand{\Ftil}{{\widetilde{F}}}
\newcommand{\Ktil}{{\widetilde{K}}}
\newcommand{\Fbar}{{\overline{F}}}
\newcommand{\fa}{{\mathfrak{a}}}
\newcommand{\bA}{{\mathbb A}}

\newcommand \bS {{\mathbb S}}

\begin{document}

\author[Cantoral Farf\'an]{Victoria Cantoral Farf\'an}
\address{Mathematische Institut, G\"ottingen \\
G\"ottingen, 37073, Germany}
\email{victoria.cantoralfarfan@mathematik.uni-goettingen.de}

\author[Li]{Wanlin Li}
\address{Department of Mathematics, Vanderbilt University\\
Nashville, TN 37240, USA}
\email{wanlin.li@vanderbilt.edu}

\author[Mantovan]{Elena Mantovan}
\address{Department of Mathematics, California Institute of Technology,
Pasadena, CA 91125, USA}
\email{mantovan@caltech.edu}

\author[Pries]{Rachel Pries}
\address{Department of Mathematics, Colorado State University\\ 
Fort Collins, CO 80523, USA}
\email{pries@colostate.edu}

\author[Tang]{Yunqing Tang}
\address{Department of Mathematics, University of California at Berkeley\\
Berkeley, CA 94720, USA}
\email{yunqing.tang@berkeley.edu}

\title{Positive density of primes of ordinary reduction
for abelian varieties of simple signature}

\maketitle

\begin{abstract}
By a result of Serre, 
if $A$ is an elliptic curve without CM defined over a number field $L$, 
then the set of primes of $L$ for which 
$A$ has ordinary reduction has density $1$.
Katz and Ogus proved the same is true when $A$ is an abelian surface, after possibly passing to a finite extension of $L$. More recently, Sawin computed the density of the set of primes of $L$ for which an abelian surface 
$A$ has ordinary reduction, depending on the endomorphism algebra of $A$.
In this paper, we prove some generalizations of these results when 
$A$ is an absolutely simple
abelian variety of arbitrary dimension whose endomorphism algebra is a CM field $F$, under specific conditions on the signature of the multiplication action of $F$ on $A$. We include explicit examples from Jacobians of curves of genus 
three through seven admitting cyclic covers to the projective line.

MSC20: primary 11G10, 11G18, 14K15; secondary 11G20, 14G35, 14H10


 








\end{abstract}

\section{Introduction}

\subsection{Serre's ordinary reduction conjecture}\label{sec_int_ord}
Let $A$ be an abelian variety defined over a number field $L$.
For a prime $\vvv$ of $L$, the reduction of $A$ modulo $\vvv$ is \emph{ordinary}
if the only slopes of its Newton polygon are $0$ and $1$.
It is conjectured that, after possibly passing to a finite extension of $L$, the set of primes for which $A$ has ordinary reduction is of density $1$; see for instance, the beginning of \S 7 in \cite{Pink98}. This conjecture is usually attributed to Serre as this question appears to be first considered by Serre in \cite[\#133]{SerrecollectedIV}; see for instance \cite[Conjecture 1.1]{Suh20}.\footnote{Throughout the paper, density means natural density of primes; some references that we cited stated their results in terms of Dirichlet density, but for the results that we cite, their proofs actually showed the results for natural density.}

This conjecture is proved when $A$ is an elliptic curve by Serre \cite[\#133]{SerrecollectedIV} or an abelian surface by an idea of Katz explained by Ogus in \cite[Proposition~2.7, Corollary~2.9]{Ogus82}. See \S\ref{sec_previous} for other known cases of this conjecture.

Moreover, as a motivation for our results, we remark that the density of ordinary reduction primes of $A$ over $L$ in the above two cases are known.
When $A$ is an elliptic curve, Serre proved that the density of ordinary primes is $1$ if the endomorphisms of $A_\Lbar$ are all defined over $L$ (see for instance \cite{Serre81}). Otherwise, the density is $1/2$; this is a direct consequence of the Shimura--Taniyama formula.
When $A$ is an abelian surface, 
Sawin \cite{Sawin16} proved that the density of the set of ordinary primes is $1$, unless $A$ has an isogeny factor with complex multiplication (CM); 
in general, the set of primes for which the abelian surface $A$ has ordinary reduction has density either $1$, $1/2$, or $1/4$. 

In this paper, we prove a special case of Serre's conjecture for certain higher dimensional abelian varieties. Let $F$ denote $\End^0(A_{\overline{L}})$, the endomorphism algebra of $A$ base changed to an algebraic closure $\overline{L}$ of $L$. We assume that\footnote{Our assumption excludes the case when $F$ is an imaginary quadratic field and $\cf=(1,1)$ because the abelian surface $A$ with $\End^0(A_{\overline{L}})\supset F$ is geometrically isogenous to the self-product of an elliptic curve and thus $\End^0(A_{\overline{L}})\neq F$.}
\begin{itemize}
    \item $F$ is a CM field, and
    \item the signature $\cf$ of $A$ is \emph{simple} in the sense of Harris--Taylor (see \cite[pages 52-53]{HarrisTaylor}).\footnote{Unlike in \cite{HarrisTaylor}, $F$ is not necessarily the composition of an imaginary quadratic field and a totally real field; such signature type is called fake Drinfeld type in \cite[Example 2.3]{RSZ}.} 
\end{itemize}
More precisely, the second assumption means that fixing an embedding $L \hookrightarrow \CC$, there exist a CM type $\Phi$ of $F$ and $\sigma_1\in \Phi$ such that for every $\sigma \in \Phi -\{\sigma_1\}$, the subspace of $\Lie (A_\CC)$ where $F$ acts via $\sigma: F\hookrightarrow \CC$ is $0$ and the subspace of $\Lie (A_\CC)$ where $F$ acts via $\sigma_1$ is of dimension $1$.

Our first result is the following.

\begin{theorem*}[Theorem \ref{Tdensityord}] \label{Thm1} 
Let $A$ be an abelian variety defined over a number field $L$.
Assume $\mathrm{End}^0(A_\Lbar) \simeq F$ is a CM field and the signature of the multiplication of $F$ on $A$ is simple.
Then, after a finite extension $L'/ L$, the set of primes of $L'$ for which the reduction of $A$ is ordinary has density $1$. Moreover, if $F/\QQ$ is Galois, then one may take $L'=LF$ and if moreover the relative dimension $2\dim A/[F:\bQ]\geq 2$, then the density of ordinary primes over $L$ is $1/[LF:L]$. 
\end{theorem*}

\subsection{The $\mu$-ordinary Newton polygon and a refinement of Serre's conjecture}
In this paper, we also provide more information on the reductions of $A$ modulo primes $\vvv$ of $L$ which are not totally split in $L'/L$ in the above theorem.
Let $p$ be the residue characteristic of $\vvv$.
The splitting behavior of $p$ in $F/\mathbb{Q}$ places constraints on the Newton polygon 
of the reduction of $A$ modulo $\vvv$ (see \cite{Kott1,Kott2} and \cite{Rapoport-Richartz}). The $\mu$-ordinary Newton polygon is the lowest Newton polygon that may occur under this restriction;\footnote{If we restrict to abelian varieties whose Hodge cycles are all generated by endomorphisms, we view $A$ as a point on a PEL type Shimura variety and the $\mu$-ordinary Newton polygon is the one that occurs in the largest Newton stratum at $p$; more generally, we consider the $\mu$-ordinary Newton polygon corresponding to the largest Newton stratum in the mod $p$ special fiber of the Hodge type Shimura variety defined by the Mumford--Tate group of $A$.} it equals the ordinary Newton polygon if $p$ is totally split in $F/\QQ$.

As a refinement of Serre's conjecture, one may ask whether the set of primes of $L$ with $\mu$-ordinary reduction is of density $1$. The proofs of Serre and Sawin for elliptic curves and abelian varieties give affirmative answers in both cases. Our second theorem gives an affirmative answer to the abelian varieties in the previous theorem under one extra assumption.

\begin{theorem*}[Theorem \ref{densitysigma}]\label{Thm2}
With notation as in the previous Theorem, assume also that 
$F/\QQ$ is Galois with abelian Galois group.
Then the set of primes of $L$ for which the reduction of $A$ is $\mu$-ordinary has density $1$.
\end{theorem*}

As motives, abelian varieties with simple signature discussed in our main theorems above are closely related to K3 surfaces (see \Cref{GUtoK3}). Independently, in an unpublished manuscript \cite{Sawin}, Sawin proved that a K3 surface over a number field has $\mu$-ordinary reduction at a density one set of primes. Sawin's method and ours are very similar.  See \Cref{rmk_Sawin} for a comparison of the results.

\subsection{Previous work}\label{sec_previous}
Here are some previous results about ordinary reduction.
In \cite[Theorems 1.7 and 2.2]{Noot95}, Noot proved that, given a group $G$ that is the Mumford--Tate group of an abelian variety over ${\mathbb C}$, there exists an abelian variety defined over a number field having Mumford--Tate group $G$ and good ordinary reduction at a set of primes of Dirichlet density $1$. 
For an abelian variety $A$ over a number field $L$ such that $\End(A_\Lbar) \simeq \ZZ$, under a condition on the Mumford--Tate group of $A$ being small, 
Pink proved that the set of primes of ordinary reduction of $A$ has Dirichlet density $1$ after a finite extension of $L$ \cite[Corollary~7.2]{Pink98}.\footnote{The condition is that the Lie algebra of the Mumford--Tate group does not have simple factors of type $C_r$ for $r\geq 3$; in particular, the Mumford--Tate group is strictly smaller than $\GSp_{2g}$ for $g\geq 3$.}
By extending the ideas of Serre, Katz and Ogus, 
Fit\'e \cite{Fite20} obtained various results on positive density of ordinary reductions
for certain abelian varieties of small dimension when $\End^0(A_\Lbar)$ contains an imaginary quadratic field with certain signature conditions; he also obtained some lower bounds on the number of Frobenius eigenvalues with $p$-adic valuation $0$ for abelian varieties whose endomorphism algebra contains a Galois field of degree $\dim A$. 
Wang and Zhang \cite{WZ25} proved that for a $\GL_2$-type abelian variety $A$ (i.e., $\End^0(A)$ contains a field of degree $\dim A$), the set of ordinary primes of $A$ is of positive density under certain assumptions on the definition field and $\End^0(A)$.
Suh \cite{Suh20} went beyond the classical method and proved the existence of infinitely many ordinary reductions for the Andr\'e motives associated to new normalized Hilbert eigencuspforms of parallel weight $(2,\dots, 2)$. Recently Hui \cite{Hui} proved that a non-CM abelian variety or K3 surface over a number field has a density $0$ set of supersingular reductions.

\subsection{Ideas of the proof}
Our proofs are based on the ideas of Serre, Katz, Ogus, and Sawin from the papers mentioned above.
The key novelty of our work lies in the proof of the second theorem, which we explain here. These authors constructed certain functions on the Galois group (arising from traces of the Frobenius action on certain cohomology groups) to detect the ordinary property; more precisely, by using the Weil bounds, at the non-ordinary points, these functions evaluate to at most finitely many values, which are independent of the prime;
then one concludes by applying the Chebotarev density theorem to each connected component of the Zariski closure of the Galois group. To obtain our second theorem, we construct functions to detect the $\mu$-ordinary property (these functions depend on the splitting behavior of $p$ in $F/\QQ$)
in Section~\ref{Sdetect}; and similarly use the Weil bound to show that 
these functions evaluate to at most finitely many values at the non-$\mu$-ordinary points.
In order to study the behavior of each connected component, Sawin used the classification of algebraic Sato--Tate groups by Fit\'e--Kedlaya--Rotger--Sutherland \cite{FKRS}. Although we do not have a full classification for the abelian varieties in our theorems, we obtain enough information on the above mentioned functions on each connected component, based on Vasiu's theorem \cite{Vasiu} that the Mumford--Tate conjecture holds for these abelian varieties and the observation that all such connected components are contained in the normalizer of the Mumford--Tate group in $\GSp_{2g}$, which in turn admits the Mumford--Tate group as a finite index subgroup. The details are in Section~\ref{Sdensity}. Section~\ref{sec2} explains the assumptions in our main theorems and discusses preliminaries for the proofs. The reader may skip most of Section~\ref{sec_assumptions} for a first reading and refer back later. Section~\ref{SECcurves} provides some applications of our theorems to Jacobians. 

We remark that it is possible to generalize these results to the case 
when $\mathrm{End}^0(A_\Lbar)$ is a division algebra over a CM field.
We leave this generalization to an interested reader.

\subsection*{Acknowledgements}
\thanks{
We thank Jeff Achter, Nir Elber, Mark Kisin, Salim Tayou, Ziquan Yang, and Rong Zhou for helpful discussions and comments. We thank Will Sawin for sharing with us his unpublished manuscript \cite{Sawin}.
Li was partially supported by NSF grant DMS-23-02511.
Mantovan was partially supported by NSF grant DMS-22-00694.
Pries was partially supported by NSF grant DMS-22-00418. 
Tang was partially supported by NSF grant DMS-22-31958 and a Sloan fellowship. 
We would also like to thank the American Institute of Mathematics for their support through the Square program.
Part of the work was done during the Spring 2023 ``Diophantine Geometry'' special semester at MSRI/SLMath; Li and Tang spent the semester there; Cantoral Farf\'an, Mantovan and Pries visited during workshops. We would like to thank AIM and MSRI/SLmath for providing stimulating working environments.
}

\section{Preliminaries and Notation}\label{sec2}

Let $A$ be an abelian variety of dimension $g$ defined over a number field $L$; let $\overline{L}$ denote an algebraic closure of $L$.
Let $\End^0(A_\Lbar)=\End(A_\Lbar)\otimes_\ZZ\QQ$ denote the (geometric) endomorphism algebra of $A$. Throughout the paper, we assume that
\begin{assumption}[F]\label{assF}
   The endomorphism algebra $\End^0(A_\Lbar)$ is isomorphic to a CM field $F$.
\end{assumption}

\begin{notation}\label{setting}
Let $*$ denote complex conjugation on $F$.
Let $F_0$ be the maximal totally real subfield of $F$. 
Let $d:=[F_0:\QQ]$.
Let $n:=2 \dim(A)/[F:\QQ] =g/d$ be the {\em relative dimension} of $A$.
\end{notation}

Fix an embedding $L\hookrightarrow\CC$. Consider the action of $F$ on ${\rm Lie}(A_\CC)$.
For each embedding $\sigma \in {\rm Hom}_\QQ(F,\CC)$,  
we denote by ${\rm Lie}(A_\CC)_\sigma$ the subspace of ${\rm Lie}(A_\CC)$ 
where $F$ acts via $\sigma:F\hookrightarrow \CC$.
We define the {\em signature} $\cf$ of $A$ (with respect to the multiplication of $F$) by
\begin{equation}\label{sig}
    \cf:{\rm Hom}_\QQ(F,\CC)\to \ZZ_{\geq 0},\quad \cf(\sigma):=\dim_\CC {\rm Lie}(A_\CC)_\sigma.\end{equation}
Note that $\cf(\sigma)+\cf(\sigma^*)=n$.

\subsection{Unitary Shimura data}\label{subsec: PEL Shimura datum}

We refer to \cite[Section 5]{Kottwitz} for the material in this subsection and the next (see also \cite[Section 2]{RSZ} for details of unitary Shimura data).

Fix a polarization $\lambda$ on $A$ also defined over $L$. 
Then the $\lambda$-Rosati involution on $\End^0(A_\Lbar)$ agrees with complex conjugation $*$ on $F$~\cite[Thm. 2, Chap. IV, \S 21 ]{Mumford74}. 
(Note that the statements are independent of the choice of a polarization defined over $L$.)

Set $V:=H^1(A_\CC,\QQ)$, equipped with the $\lambda$-Weil pairing $\langle \cdot, \cdot\rangle$, which is a symplectic form over the $\QQ$-vector space $V$. Due to the $F$-action, we may view $V$ as an $F$-vector space and it is endowed with a Hermitian form $(\cdot, \cdot )$ such that $\langle \cdot, \cdot\rangle=\tr_{F/\QQ} \xi (\cdot, \cdot )$, where  $\xi\in F$ is a non-zero totally imaginary element.
Consider the Hodge structure $h:\CC\to \End_{F\otimes \RR}(V_\RR)$. 
The tuple $\mathcal{D}=(F, *, V, (\cdot,\cdot), h)$ defines a {\em Shimura datum of PEL-type}, associated with the algebraic group $G$ over $\QQ$, the unitary similitude group with similitude in $\mathbb{G}_m$ defined as follows: for every $\QQ$-algebra $R$, we have
\begin{equation} \label{Edefgroup}
G(R):=\{x\in \GL_F(V\otimes_\QQ R) \mid \exists c(x)\in R^\times, \text{such that } (xv, xw)=c(x)(v,w) \forall v, w \in V\otimes R \}.
\end{equation}

Let $E \subset \CC$ denote the {\em reflex field} of the Shimura datum.  By definition, $E$ is a subfield of the Galois closure of $F/\QQ$. More precisely, $E$ is the fixed field of the subgroup of $\Gal(\overline{\QQ}/\QQ)$ fixing the signature type $\cf$; that is $\Gal(\overline{\QQ}/E)=\{\gamma\in \Gal(\overline{\QQ}/\QQ) \mid \cf(\gamma\circ \sigma)=\cf(\sigma),\forall \sigma:F\hookrightarrow \overline{\QQ}\}$ (here as all images of $\sigma:F \hookrightarrow \CC$ lie in $\overline{\QQ}\subset \CC$, we then consider $\sigma: F \hookrightarrow \QQ$).

Let $G_0$ denote the unitary group $U(V,(\cdot,\cdot))$ over $F_0$ and let $G_1:=\Res_{F_0/\QQ}G_0\subset G$.
At each place $F_0\hookrightarrow \RR$, the unitary group $G_0$ has signature equal to $(\cf(\sigma),\cf(\sigma^*))$, for $\sigma:F\hookrightarrow \CC$ satisfying $\sigma_{\vert F_0}$ is the given embedding $F_0\hookrightarrow \RR$ and $\sigma(\xi^{-1})$ has positive imaginary part.   

Let $\adeles_f$ be the ring of finite ad\`eles of $\QQ$. 
For an open compact subgroup $\mathcal{U}\subset G(\adeles_f)$, 
consider the associated moduli space $S_{\mathcal{U}}=S_{\mathcal{U}}(\mathcal{D})$
of polarized abelian varieties of level $\mathcal{U}$.
Then $S_{\mathcal{U}}$ is a smooth Deligne--Mumford stack defined over $E$.

\begin{definition} \label{pgood}
We call a rational prime $p$ {\em good } if $p$ is odd and: \begin{itemize}
\item $A/L$ has good reduction at $\vvv$ for all primes $\vvv$ of $L$ above $p$; 
\item the degree of the polarization $\lambda$ is relatively prime to $p$; and 
\item $p$ is unramified in the extension $F/\QQ$. 
\end{itemize}

We call a prime $\vvv$ of $L$ {\em good} if $\vvv$ lies above a good rational prime $p$. 
\end{definition}
A good prime $p$ is unramified in the reflex field $E/\QQ$
and so is a prime of good reduction for the Shimura datum $\mathcal{D}$.
Furthermore, if $\mathcal{U}=\mathcal{U}^p\mathcal{U}_p$ where 
$\mathcal{U}^p\subset G(\adeles_f^{(p)})$ is open compact and $\mathcal{U}_p\subset G(\QQ_p)$ is 
hyperspecial, then $S_{\mathcal{U}}$ admits a canonical smooth model over $\mathcal{O}_E\otimes_\ZZ \ZZ_{(p)}$.

\subsection{The $\mu$-ordinary Newton polygon} \label{SmuordNP}
Let $\vvv | p$ be a good prime of $L$ for $A$. 
Let $\FF$ be the residue field of $\vvv$, which is a finite field of characteristic $p$; let $m$ be the degree of $\FF$ over $\FF_p$.
Let $W(\FF)$ denote the Witt vector ring of $\FF$. 

Let $A_\FF$ denote the reduction of $A$ modulo $\vvv$. 
Let $\varphi$ denote the Frobenius action on the crystalline cohomology group 
$H^1_{\rm cris}(A_\FF/W(\FF))$. 
The {\em Newton polygon} $\nu(A_\FF)$ of $A_\FF$ is the multi-set of rational numbers $\lambda$, called the {\em slopes},
such that $m\lambda$ are the $p$-adic valuations of the eigenvalues of the $W(\FF)$-linear map $\varphi^m$ (here we normalized the valuation $v$ such that $v(p)=1$).

The \emph{ordinary} Newton polygon is the symmetric Newton polygon with slopes $0$ and $1$.
Given two Newton polygons $\nu_1$ and $\nu_2$, we write $\nu_1\leq \nu_2$ (resp. $\nu_1<\nu_2$) if the lower convex hull of slopes $\nu_1$ is on or above that of $\nu_2$ (resp. $\nu_1$ is strictly above $\nu_2$), and the two polygons have the same endpoints.  

Under Assumption F (\ref{assF}), 
in \cite[Theorem 4.2]{Rapoport-Richartz}, Rapoport and Richartz described the lowest polygon in terms of the splitting behavior of $p$ in $F/\QQ$; in \cite[Theorem 1]{wedhorn}, Wedhorn proved that there exists an abelian variety defined over a number field with multiplication by $F$ and good reduction at $\vvv$ whose reduction admits the lowest polygon in \cite{Rapoport-Richartz} as its Newton polygon; in other words, this lowest polygon describes the dense open Newton stratum on the special fiber of $S_{\mathcal{U}}$.
This lowest polygon is called the {\em $\mu$-ordinary} Newton polygon; it equals the ordinary Newton polygon if $p$ is totally split in $F/\QQ$.

\subsection{The Mumford--Tate group and the $\ell$-adic monodromy group}\label{sec_MTA}

Let $M_A$ denote the Mumford--Tate group of $A$; since it fixes all Hodge cycles, we have $M_A\subseteq G$, with $G$ defined in \eqref{Edefgroup}.

Fix a prime $\ell$.  Let
$V_\ell := H^1_{\text{\'et}}(A_{\bar{L}},\QQ_\ell)$ and let $\rho_{A,\ell}$ denote the Galois representation
\begin{equation} \label{Erho}
\rho_{A,\ell}: \Gal(\overline{L}/L)\to \GL(V_\ell). 
\end{equation}
The \emph{ $\ell$-adic monodromy group $\Gl$} of $A$ is the Zariski closure of the image of $\rho_{A,\ell}$. 
Let $\Glcon$ denote the identity component of $\Gl$, and $\pi_0(\Gl):=\Gl/\Glcon$ the component group of $\Gl$. 
Let $\chi$ denote the symplectic similitude character of $\GSp(V)$ (and hence of $\GSp(V_\ell)$) and we also use $\chi$ to denote the induced similitude characters of $G$ and $G_{A,\ell}$. We define $G_{A,\ell}^{\chi=1}:= \Gl \cap \Sp(V_\ell)$.
The Mumford--Tate Conjecture predicts that $\Glcon=M_{A,\QQ_\ell}$.  

\begin{definition}\label{DefLconn}
Let $L^{\rm conn}/L$ be the field extension of $L$ such that $\rho_{A,\ell}^{-1}(\Glcon)=\Gal(\overline{L}/L^{\rm conn})$. 
\end{definition}
By construction, $L^{\rm conn}/L$ is a finite Galois extension and $\Gal(L^{\rm conn}/L)$ is canonically isomorphic to $\pi_0(\Gl)$.
By the work of Serre \cite[\#133,p.15, Th\'eor\`eme, \#135, 2.2.3]{SerrecollectedIV}(see also \cite[Theorem~0.1]{LarsenPink}), the extension $L^{\rm conn}/L$ is independent of $\ell$. The reader may also see \cite[\S 5]{Pink98} for standard results on $M_A$ and $\Gl$.

\subsection{Running assumptions and basic properties}\label{sec_assumptions}
We list various assumptions on $A$ used in our main results; recall that $A$ satisfies Assumption F.  

\begin{assumption}[S]\label{assS} 
   The signature $\cf$ is simple (as in \cite[Lemma I.7.1, page 52]{HarrisTaylor});
explicitly, 
there exists a CM type $\Phi$ of $F$ and $\sigma_1\in \Phi$ such that 
\begin{equation} \label{Easspart2}
  \cf(\sigma_1)=1,\text{ and } \cf(\sigma)=0\text{ for all }\sigma\in\Phi-\{\sigma_1\};
\end{equation}
\end{assumption}
Note that geometrically simple CM abelian varieties satisfy Assumptions F and S. 
\begin{assumption}[A]\label{assA}  
    The extension $F/\QQ$ is Galois with abelian Galois group.
\end{assumption}

\begin{assumption}[D]\label{assD}  
 For $F/\bQ$ abelian,   \label{assLconn} $L^{\rm conn} \subset FL$, where $L^\conn$ is defined in \Cref{DefLconn}.
\end{assumption}
When $F/\QQ$ is not abelian, we replace the above Assumption D by the following: let $E$ denote the reflex field of the Hodge type Shimura variety associated to $M_A$; note that $E$ is also the reflex field of the PEL Shimura variety defined by the unitary similitude group $G$. 
We have $ E \hookrightarrow \overline{L}\hookrightarrow \CC$, where $\overline{L}$ denotes the algebraic closure of $L$ in $\CC$ under the fixed embedding $L \hookrightarrow \CC$. 
The corresponding Assumption D is that $L^{\rm conn} \subset EL$ (the composite field is taken with respect to the above embedding of $E$ into $\overline{L}$).

Note that $E$ is always contained in the Galois closure of $F$ over $\bQ$. We now assume the extra Assumption S. When $n\geq 2$, via $\sigma_1$ in Assumption S, $F\hookrightarrow E$. 
When $F/\bQ$ is Galois and $n\geq 2$, or when $F/\bQ$ is abelian, we have $F=E$.
Note that $E$ may not equal $F$ when $F$ is not Galois over $\bQ$.\\

Consider the natural group homomorphism $G \rightarrow \Res_{F/\bQ} \bG_m$ given by $G\subset \Res_{F/\bQ} \GL_n \xrightarrow{\det} \Res_{F/\bQ} \bG_m$. Let $T^F$ be the algebraic subgroup of $\Res_{F/\bQ} \bG_m$ defined by 
\[T^F(R):=\{x\in (R\otimes_{\bQ} F)^\times \mid \overline{x} \cdot x \in R\}\] 
for every $\bQ$-algebra $R$, where $\overline{\cdot}$ denotes the involution on $\Res_{F/\bQ} \bG_m$ induced by complex conjugation on $F$. By definition of $G$, the map $G \rightarrow \Res_{F/\bQ} \bG_m$ factors through and surjects onto $T^F$.
\begin{assumption}[C]\label{assC}
The induced map $\pi: M_A \hookrightarrow G \rightarrow T^F$ is surjective.
\end{assumption}
We provide several equivalent versions of Assumption C. Note that the natural inclusion $\bG_m \subset \GL_n$ given by $a \mapsto \diag (a, \dots, a)$ induces inclusions $\Res_{F/\bQ} \bG_m \subset \Res_{F/\bQ} \GL_n$ and $T^F \subset G$. Note that the center $Z(G)=T^F$.
Assumption F and Hodge theory imply that the center $Z(M_A)\subset T^F$ and thus Assumption C is equivalent to $Z(M_A)=Z(G)$.

Another equivalent condition is described in terms of the signature $\cf$. Let $X_*(T)$ denote the space of cocharacters of an algebraic torus $T$ (base changed to an algebraically closed field). We then have a natural identification $X_*(\Res_{F/\bQ} \bG_m)\otimes \bQ = \oplus_{\sigma \in \Hom_{\bQ}(F,\Qbar)} \bQ_{\sigma}$ (here $\bQ_\sigma$ denotes the $1$-dimensional $\bQ$-vector space associated to $\sigma: F \hookrightarrow \Qbar$). Under this identification, $X_*(T^F)\otimes \bQ \subset X_*(\Res_{F/\bQ} \bG_m)\otimes \bQ$ is the subspace \[\{(x_{\sigma}) \mid x_\sigma + x_{\sigma^*} \text{ is a constant independent of }\sigma \}.\] Every signature $\cf$ can be viewed as an element $(f(\sigma))_{\sigma \in \Hom_{\bQ}(F,\Qbar)} \in X_*(T^F)\otimes \bQ$ and we still use $\cf$ to denote this element (here we pick an embedding $\Qbar \hookrightarrow\CC$ to identify $\Hom_{\bQ}(F,\Qbar)=\Hom_{\bQ}(F,\CC)$; the formulation of the final conclusion is independent of the choice of this embedding). The $\Gal(\Qbar/\bQ)$-action on $X_*(T^F)\otimes \bQ \subset X_*(\Res_{F/\bQ} \bG_m)\otimes \bQ$ is given by $(y_\sigma):=\tau. (x_\sigma)$ with $y_{\tau\circ\sigma}= x_\sigma$.

\begin{lemma}\label{AssCeqv}
    Assumption C is equivalent to the condition $ \spn_{\bQ} \{\tau \cf \mid \tau \in \Gal(\Qbar/\bQ)\}= X_*(T^F)\otimes \bQ$; in other words, it is equivalent to the condition that $ span \{\tau \cf \mid \tau \in \Gal(\Qbar/\bQ)\}$ is $([F_0:\bQ]+1)$-dimensional.
\end{lemma}
\begin{proof}
   Let $\mu:\bG_{m,\CC} \rightarrow M_{A,\CC}$ denote the Hodge cocharacter of $A_{\CC}$. By definition of $\cf$, we have $\pi\circ \mu \in X_*(T^F)$ equals $\cf$. By definition, the Mumford--Tate group $M_A$ is generated by all $\Gal(\Qbar/\bQ)$-conjugates of all $M_A(\CC)$-conjugates of $\mu$. (More precisely, all $M_A(\CC)$-conjugacy classes of cocharacters are defined over $\Qbar$.  Hence we may pick a cocharacter $\mu$ over $\Qbar$ in the same conjugacy class as the Hodge cocharacter, and then $M_A$ is generated by the $M_A(\Qbar) \rtimes \Gal(\Qbar/\bQ)$-conjugates of $\mu$.) Since $T^F$ is commutative, the image $\pi(M_A)$ is then generated by all $\Gal(\Qbar/\bQ)$-conjugates of $\pi\circ\mu$, which gives the desired equivalence.
\end{proof}

In the degenerate case $n=1$ (i.e., when $A$ is a geometrically simple CM abelian variety), the Mumford--Tate group $M_A$ is an algebraic torus whose rank is between $[F_0:\QQ]+1$ and $2+\log_2 [F_0:\QQ]$ (proved by Ribet \cite[\S 3]{Ribet}); moreover, the work of Ribet \cite[\S 3]{Ribet} and \cite{Ribet2} (see also Dodson \cite[\S 1]{Dodson} for a summary) gives various sufficient conditions for $M_A=G$ and also examples when $M_A\neq G$.

Under Assumptions F and S, there are examples beyond the CM case when $Z(M_A)\neq Z(G)$ (i.e., Assumption C fails). 
\begin{example}
    Suppose $F=\QQ(\zeta_7)$ and the signature is $\cf=(1,0,0,2,2,1)$. (Note that all embeddings $F\hookrightarrow \CC$ are given by $\sigma_i: \zeta_7 \mapsto \zeta_7^i$ for $i=1,\dots, 6$ and we write $\cf$ as a vector $(\cf(\sigma_1),\dots, \cf(\sigma_6))$).  Then $n=2$ and we have
\[\spn_{\bQ} \{\tau \cf \mid \tau \in \Gal(\Qbar/\bQ)\}=\spn_{\bQ}\{(1,0,0,2,2,1), (2,1,2,0,1,0), (2,0,1,1,2,0), (1,\dots, 1)\},\]
where the first three vectors are $\tau\cf$ for $\tau:\zeta_7 \rightarrow \zeta_7^i$, where $i=1,2,3$ and $(1,\dots,1)$ is included because $\cf + c\cf=(n,\dots, n)$, where $c$ denotes complex conjugation. We note that this is a $3$-dimensional $\bQ$-vector space and hence by \Cref{AssCeqv}, $Z(M_A)\neq Z(G)$. In particular, $M_A\neq G$, there are extra Hodge cycles beyond the ones generated by endomorphisms. 
\end{example}
On the other hand, by \Cref{AssCeqv}, Assumption C holds for $F=\QQ(\zeta_7)$ and $\cf=(1,0,0, 3,3,2)$ and we have $M_A=G$ by the following lemma.

\begin{lemma}\label{MderisGder}
Assumptions F and S imply that $M_A^{\der}=G^{\der}$, where $G$ is defined in \eqref{Edefgroup}. As a consequence, under Assumptions F and S, we have that Assumption C is equivalent to $M_A=G$.
\end{lemma}

\begin{proof}
The work of Serre \cite[\S 3]{Serre79} and Pink \cite[\S\S 4-5]{Pink98} provides a description of connected reductive groups which have certain structure similar to the Mumford--Tate group of an abelian variety. We follow the description and ideas in Pink's paper to prove the lemma by showing that Assumptions F and S force the derived part of any possible Mumford--Tate group to be $G^{\der}$.\footnote{The proof of this lemma is independent of the rest of the paper, thus we do not recall Pink's paper in full detail. In terms of the geometry of Shimura varieties, the condition $M_A^{\der}=G^{\der}$ is equivalent to that the corresponding point of $A$ (after choosing suitable level structure) on the PEL type Shimura variety $S_{\cU}$ does not lie on any Shimura sub-variety, i.e., $A$ is Hodge generic. The reader may skip this lemma and replace Assumption F by the assumption that $A$ is Hodge generic; this does not affect the rest of the paper.}

Let $\Ftil\subset \Fbar$ denote the Galois closure of $F$ over $\bQ$ in $\Fbar$. Let $\rho$ denote the standard representation of $M_A, M_A^{\der}$ and $G, G^{\der}$ on $V$; let $\rho_{\Ftil}$ denote the base change of $\rho$ to $\Ftil$. By definition of $G$, we have $\rho_{\Ftil}=\oplus_{\sigma \in \Hom_{\bQ}(F,\Ftil)} \rho_{\sigma}$, where each $\rho_{\sigma}$ is an absolutely irreducible representation of $G_{\Ftil}$ of dimension $n$ and the image $\rho_{\sigma}(G_{\Ftil})=\GL_{n,\Ftil}$. Therefore $\rho_\sigma(G_\Ftil^{\der})=\SL_{n,\Ftil}$ and $\rho_\sigma(Z(G)_\Ftil)=\bG_{m,\Ftil}$. We then conclude that $\rho_\sigma$ is an absolute irreducible representation of $M^{\der}_{A,\Ftil}$. Indeed, if not, the centralizer of $\rho_{\sigma}(M^{\der}_{A,\Ftil})$ in $\GL_{n,\Ftil}$ would contain a torus of rank at least $2$; on the other hand, since $\rho_\sigma(Z(M_A)_\Ftil) \subset \rho_\sigma(Z(G)_\Ftil)=\bG_{m,\Ftil}$, we have that the centralizer of $\rho_\sigma(M_{A,\Ftil})$ is the same as that of $\rho_{\sigma}(M^{\der}_{A,\Ftil})$, which then contradicts that the centralizer of $M_A$ in $\End(V)=\End(A_\Lbar)^0=F$. 

In the sense of \cite[Definition 4.1]{Pink98}, $(M_A, \rho)$ is a strong Mumford--Tate pair of weights $\{0,1\}$ over $\bQ$ with respect to the Hodge cocharacter $\mu$ (see also \cite[Fact 5.9]{Pink98}; note that as in the proof of \Cref{AssCeqv}, we pick $\mu$ defined over $\Qbar$ in the same conjugacy class as the Hodge cocharacter directly arising from $A_\CC$). Therefore, with respect to the projection to $\rho_\sigma (M_A)$ of all $\Gal(\Qbar/\bQ)$-conjugates of $\mu$, we have that $(\rho_\sigma(M_A), \rho_\sigma)$ is a weak Mumford--Tate pair of weights $\{0,1\}$. We can also check directly that with respect to the same cocharacters, $(\rho_\sigma(G), \rho_\sigma)$ is a weak Mumford--Tate pair of weights $\{0,1\}$. 

Recall that $\rho_\sigma$ is an absolute irreducible representation of $\rho_\sigma(M_A)$. Therefore by the discussion in \cite[pp.~210-211]{Pink98}, $\rho_\sigma(M_A)_\Fbar=\bG_{m,\Fbar}M_1\cdots M_r$ and $\rho_\sigma \simeq \rho_0 \boxtimes \cdots \boxtimes \rho_r$, where $M_1,\dots, M_r$ are pairwise distinct (almost) simple factors of $\rho_\sigma(M_A)^{\der}_\Fbar$ and the above product is an almost direct product, and moreover $\rho_0$ is the representation of $\bG_{m,\Fbar}$ as scalar matrices and $\rho_i$ is an irreducible representation of $M_i$. Furthermore, for each $i$, Pink shows that $(\bG_{m,\Fbar}M_i, \rho_0\boxtimes \rho_i)$ is a (weak) Mumford--Tate pair of weights $\{0,1\}$ with respect to some $\Gal(\Qbar/\bQ)$-conjugates of $\mu$. More precisely, Pink shows that at least one of these cocharacters, denoted by $\mu_i$, has non-trivial component $M_i$ and moreover, this $\mu_i$ must be trivial in the other $M_j$'s for $j\neq i$; therefore we can naturally view $\mu_i$ as a cocharacter of $\bG_{m,\Fbar}M_i$. Over an algebraically closed field (of characteristic $0$), the classification of (weak) Mumford--Tate pairs of weights $\{0,1\}$ such that the derived groups are almost simple is given in Serre \cite[\S 3, Annexe]{Serre79}; see also \cite[Table 4.2]{Pink98}; the classification also includes information on possible cocharacters and the multiplicities of weights $\{0,1\}$ for each case. 

Since the multiplicities correspond to signatures $\cf$, Assumption S implies that the projections of all $\Gal(\Qbar/\bQ)$-conjugates of $\mu$ to $\rho_\sigma (M_A)$ are either trivial on $\rho_\sigma (M_A)^{\der}$ or have multiplicities of weights to be $(1, n-1)$ or $(n-1,1)$. Note that if $r\geq 2$, then the multiplicities of weights of $(\rho_\sigma(M_A), \rho_\sigma)$ must be of form $(ad,bd)$, where $(a,b)$ are the multiplicities of $(\bG_{m,\Fbar}M_i, \rho_0\boxtimes \rho_i)$ and $d=\prod_{j\neq i, 1\leq j \leq r} \dim \rho_j$. Since $d\geq 2$, this contradicts the assumption that the multiplicities are $(1, n-1)$ or $(n-1,1)$; in other words, we have proved that $\rho_\sigma(M_A)^{\der}_\Fbar$ is almost simple. Again with our signature/multiplicity restriction, we compare with \cite[\S 3, Annexe]{Serre79} or \cite[Table~4.2]{Pink98} to conclude that $\rho_\sigma(M_A)^{\der}_\Fbar$ is of type A and $\rho_\sigma$ is either the standard representation or its composition with the outer automorphism (inverse transpose). Since $M_A \subseteq G$, we then conclude that $\rho_\sigma(M_A^{\der})=\rho_\sigma(M_A)^{\der}=\rho_\sigma(G^{\der})$.

Note that by definition $G^{\der}_{\Ftil^+}=\prod_{\sigma \in \iota^{-1}(\Phi)} (\rho_\sigma \oplus \rho_{\sigma^*})(G^{\der}_{\Ftil^+})$ (recall that $\Phi$ is the CM type in assumption S and we fix an embedding $\iota: \Ftil \hookrightarrow \CC$ to identify $\Hom_{\bQ}(F,\Ftil)= \Hom_{\bQ}(F,\CC)$). We have proved that $M_{A,\Ftil^+}^{\der} \subset G^{\der}_{\Ftil^+}$ maps surjectively onto each factor $\rho_\sigma(G^{\der})$ and we will use Goursat's lemma to conclude that $M_A^{\der}=G^{\der}$. More precisely, $(\rho_\sigma \oplus \rho_{\sigma^*})(G^{\der}_{\Ftil^+})$ is, up to center, a unitary group with respect to a Hermitian form of signature $(\cf(\iota\circ \sigma),\cf(\iota\circ \sigma))$. For any $\sigma, \sigma'$ such that $\sigma \neq \sigma', (\sigma')^*$, we can choose a suitable $\iota$ such that $\iota \circ \sigma = \sigma_1 \in \Phi$; then we have $(\rho_\sigma \oplus \rho_{\sigma^*})(G^{\der}_{\Ftil^+})$ and $(\rho_{\sigma'} \oplus \rho_{(\sigma')^*})(G^{\der}_{\Ftil^+})$ are not isomorphic as $\RR$-groups and thus by Goursat's lemma, $M_{A,\Ftil^+}$ surjects onto $(\rho_\sigma \oplus \rho_{\sigma^*}\oplus \rho_{\sigma'} \oplus \rho_{(\sigma')^*})(G^{\der}_{\Ftil^+})$. We then conclude that $M_A^{\der}=G^{\der}$.

We then have $M_A=Z(M_A)M_A^{\der} = Z(M_A) G^{\der} \subset Z(G) G^{\der}=G$ and by the discussion after Assumption C, we conclude that $M_A=G$ if and only if Assumption C holds.
\end{proof}

\begin{remark}\label{GUtoK3}
Suppose the point of $S_U$ represented by $A$ (over a number field) lies in a Shimura subvariety which is strictly smaller than $S_U$.  Then by \Cref{MderisGder} , $\End^0(A_{\overline{\QQ}})$ is strictly larger than $F$. Indeed, any Shimura subvariety of $S_{\cU}$ is cut out using extra endomorphisms (this can be proved using an inductive argument on dimension and \Cref{MderisGder}). Alternatively, using the work of Hofmann \cite[\S 4]{Hofmann} (see also \cite[\S\S 2.1,6.2]{BHKRY}), one obtains a morphism from each connected component of $S_{\cU}$ to some connected component of the orthogonal Shimura variety associated to the $F_0$-quadratic space $(V, \tr_{F/F_0} (\cdot, \cdot ))$;\footnote{More precisely, to obtain the weight $0$ Hodge structure parametrized by the orthogonal Shimura variety, we indeed consider $\Hom_F(V, V_0)$, where $V_0$ denotes $H^1(A_0,\QQ)$ of an abelian variety $A_0$ with complex multiplication by $(F,\Phi^*)$ (all such $A_0$ are isogenous over $\CC$, hence $V_0$ is well defined). With the induced Hodge structure, $\Hom_F(V, V_0)$ can be viewed as the Hodge structure associated to a point on the orthogonal Shimura variety.\label{CMA0}} one then obtains the above assertion on Shimura subvarieties of $S_{\cU}$ from the classification of Shimura subvarieties of the orthogonal Shimura variety, which is a direct consequence Zarhin's theorem on Mumford--Tate groups of K3 type Hodge structures \cite[Theorems 2.2.1, 2.3.1]{Zarhin}. See also Jiang \cite[\S 6, Proposition 6.10, Remark 6.11, \S\S 6.3-6.4]{Jiang}, in which he explains that all Shimura subvarieties of a GSpin Shimura variety are cut out by endomorphisms of the corresponding family of Kuga--Satake abelian varieties.
\end{remark}

\begin{remark}\label{vasiu}
By \Cref{MderisGder} and Vasiu's theorem \cite[Theorem~1.3.4 (a)]{Vasiu}, Assumptions F and S imply that the Mumford--Tate Conjecture holds for $A$.
If we further assume Assumption C, then $G_{A,\ell}^\circ=M_{A,\QQ_\ell}=G_{\QQ_\ell}$.
\end{remark}

\begin{lemma}
   Given a CM field $F$ and a CM type $\Phi$, there exists a finite set $N\subset \ZZ_{>0}$ such that all abelian varieties which satisfy Assumptions F and S with respect to the given $(F,\Phi)$ also satisfy Assumption C as long as $\dim A \notin N$.
\end{lemma}

\begin{proof}
It suffices to prove that there are at most finitely many values of the relative dimension $n$ such that Assumption C is not satisfied. We use \Cref{AssCeqv}. We observe that the dimension of $W(\cf):=\spn_{\bQ} \{\tau \cf \mid \tau \in \Gal(\Qbar/\bQ)\}$ is the same as $\dim W(2\cf)$; therefore Assumption C is equivalent to the equality $\dim W(2\cf - n)=[F_0:\QQ]$. Since $2\cf(\sigma)-n = - (2\cf(\sigma^*)-n)$, Assumption C is equivalent to the $[F_0:\QQ]\times [F_0:\QQ]$-matrix $B(n):=(\tau(2\cf)(\sigma) -n)_{\tau, \sigma \in \Phi}$ being an invertible matrix. Note that $\det B(n)$ is a polynomial in $n$.
Thus we only need to show $\det B(n)$ is not the zero polynomial, which implies that it has at most finitely many roots in $\ZZ_{>0}$. To see $\det B\not\equiv 0$, we may substitute $n=0$ and we observe that $\det B(0)=\pm 2^{[F_0:\QQ]}$ because by definition (since $\Phi$ is a CM type), each row and column of $B(0)$ has exactly one nonzero element, which is $\pm 2$.
\end{proof}

For the examples considered in \Cref{CorMoo}, all the simple factors of $J_C$ (except the ones in Moonen's family M[19]) satisfy Assumptions F, S, A, and C; and we will see that they automatically satisfy Assumption D by \Cref{ACimpliesAD}. 

\begin{lemma}\label{ACimpliesAD}
  Under Assumptions F and S, the definition field of all endomorphisms of $A_\Lbar$ is $EL$; in particular, $EL\subset L^\conn$. Therefore Assumptions F, S, and C imply that $L^{\rm conn} = EL$; i.e., Assumptions F, S and C imply Assumption D when $F/\QQ$ is abelian.\footnote{The main part of this lemma that all endomorphisms are defined over $EL$ was first pointed to us by Salim Tayou and he gave a proof using Shimura varieties as follows. As we need to pin down the level structure, we assume that $A$ is principally polarized and $\cO_F=\End(A_\Lbar)$ for simplicity (other cases shall admit a similar argument). Let $\fa$ denote the $\cO_F$-lattice $H^1(A_\CC, \ZZ)\subset V$ and set $\cU\subset G(\bA_f)$ to be the stabilizer of $\fa$. Then $S_\cU$ parametrizes principally polarized abelian varieties $B$ along with $\iota: \cO_F \hookrightarrow \End(B)$ satisfying the signature condition. Let ${\mathcal A}_g$ denote the moduli stack of principally polarized abelian varieties of dimension $g=\dim A$. By Assumptions F and S, the forgetful map $S_\cU \rightarrow {\mathcal A}_g$ is injective on the sub-locus of points having the property that $\End(B)=\cO_F$. Since this morphism of Shimura varieties is defined over $E$, we conclude by the moduli interpretation of $S_\cU$ that all endomorphisms of $A$ are defined over $EL$.}
\end{lemma}
\begin{proof}
We first prove that all endomophisms of $A_\Lbar$ are defined over $EL$.
Since $L^{\rm conn}/L$ is Galois, we have $EL^{\rm conn}/EL$ Galois and we use Galois descent to argue that all endomorphisms of $A$ are defined over $EL$. By definition of $L^\conn$, we have $F=\End^0(A_{\Lbar})=\End^0(A_{L^\conn})$ and therefore, each $\gamma \in \Gal(EL^\conn/EL)$ induces an automorphism $\gamma'$ on $F$. In other words, for any $s\in \End(A_{L^\conn})\subset F$, we have that for all $x\in A(\Qbar)$, $\gamma (s(\gamma^{-1}(x)))=(\gamma'(s))(x)$. 

Now we consider the signature of various $F$-actions on $\Lie A_\CC$. Since $\gamma$ fixes the reflex field $E$, we have that the $F$-action on $\Lie A_\CC$ via $\gamma s \gamma^{-1}$ for $s\in F$ still has signature $\cf$ as the natural $F$-action without the $\gamma$-twist. On the other hand, the $F$-action on $\Lie A_\CC$ via $\gamma' s$ for $s\in F$ has signature $\cf'$ defined as $\cf'(\sigma)=\cf(\sigma \circ \gamma')$ for all $\sigma \in \Hom_\QQ(F,\CC)$. By definition of $\gamma'$, we then have $\cf(\sigma)=\cf'(\sigma)=\cf(\sigma \circ \gamma')$. We then use Assumption S to conclude that $\sigma'=\id$. Indeed, if $n\geq 3$, since $\sigma_1$ is the unique element such that $\cf(\sigma_1)=1$, we have $\sigma_1=\sigma_1 \circ \gamma'$ and thus $\gamma'=\id$. If $n=2$, then $[F_0:\QQ]\geq 2$ by Assumption F and $\sigma_1, \sigma_1^*$ are exactly the elements with $\cf(\sigma_1)=\cf(\sigma_1^*)=1$, which implies $\gamma'=\id$ or complex conjugation. Since $[F_0:\QQ]\geq 2$, complex conjugation does not preserve $\cf$ for other embeddings and thus $\gamma'=\id$. If $n=1$, Assumption F implies that the CM type of $A$ is primitive and thus $\gamma'=\id$.

On the other hand, let $L' \subset \Lbar$ denote the definition field of all endomorphisms of $A_\Lbar$. Consider $\gamma \in \Gal(\Lbar/L')$. Note that the signature $\cf'$ associated to the $\gamma$-twisted $F$-action on $A_\CC$ above is given by $\cf'(\gamma \circ \sigma)=\cf(\sigma)$. By definition, $\gamma$ induces the identity automorphism on $F$. Thus $\cf'=\cf$ and we have $\cf(\gamma \circ \sigma)=\cf(\sigma)$ for all $\sigma$. Thus by the definition of the reflex field, $\gamma \in \Gal(\Lbar/E)$ and then we conclude that $E \subset L'$. Combined with the previous discussion, we have $L'=EL$.

The second assertion follows from the first one since by Assumption C and \Cref{vasiu}, $L^{\rm conn}$ is the minimal field extension of $L$ where all endomorphisms of $A_\Lbar$ are defined. 
\end{proof}

\begin{remark}
\Cref{ACimpliesAD} is a special property for certain PEL type unitary Shimura varieties. For more general abelian varieties $A$ over $L$, we may ask whether all geometric endomorphisms are defined over $EL$, where $E$ is the reflex field of the corresponding PEL type Shimura variety. 
This does not hold in general.
For instance, in \cite[p. 1422]{FKRS}, the authors study an abelian surface $A$ (the Jacobian of the curve $y^2=x^6+x^2+1$) defined over $\QQ$, which is isogenous over $\QQ(i)$ to a product of two elliptic curves, each without CM, and which are non-isogenous over $\overline{\QQ}$.  In this particular case, $L=\QQ$ and $E=\QQ$, but $L^{\rm conn}=\QQ(i)$.
\end{remark}

\begin{remark}\label{A45}
Although $L^{\conn}$ contains the definition field of all geometric endomorphisms, it may be strictly larger due to the existence of exceptional Hodge cycles (namely the Hodge cycles not coming from endomorphisms). 
See \cite[Theorem~1.2.1]{CLV} for examples of type III abelian varieties. We also remark that Assumptions F and S are not enough to deduce $L^{\conn}\subset EL$; see \cite[Example 4.1]{SZ98} for an absolutely simple CM abelian variety violating $L^{\rm conn} \subset FL=EL$.

On the other hand, it may still happen that although there are extra Hodge cycles, $L^{\conn}$ is still the definition field of all geometric endomorphisms. The work of Gallese--Goodson--Lombardo \cite{GGL}\footnote{We thank Nir Elber for pointing out this reference to us.} provides some concrete examples of absolutely simple CM abelian varieties over $L$ with exceptional Hodge cycles but with $L^{\conn}=EL$. For example, following the notation in \cite[\S 3, Theorem 3.0.1]{GGL}, there is a $\QQ$-simple factor $X_{21}$ of the Jacobian of the hyperelliptic curve $y^2=x^{21}+1$ with the property that $X_{21}$ is an absolutely simple CM abelian variety of dimension $6$ with CM by $\QQ(\zeta_{21})$ and the CM type is given in \cite[Lemma 3.1.3]{GGL}. We have that the Mumford--Tate group of $X_{21}$ is a rank $6$ torus and thus it admits exceptional Hodge cycles. The reference \cite[\S 4]{GGL} provides a concrete recipe to turn defining equations of the Mumford--Tate group into exceptional Hodge and Tate classes. The equation that we use here is $x_1x_4x_{10}(x_2x_5x_8)^{-1}=1$ (here we use the same notation as \cite{GGL} and $x_i$ denotes the coordinate of $(T^F)_{\Fbar}$ corresponding to the complex embedding $\zeta_{21} \mapsto \zeta_{21}^i$). The reference \cite[Theorem 6.2.2, Proposition 6.4.3]{GGL} provides an explicit formula for the definition field of these Hodge/Tate cycles in terms of Gamma functions and a direct computation using this formula shows that $L^{\conn}=\QQ(\zeta_{21})=EL=FL$. 
\end{remark}

\section{Detecting $\mu$-ordinary reduction} \label{Sdetect}

As in Section~\ref{sec2}, 
let $A$ be an abelian variety defined over a number field $L$. 
Suppose Assumptions~F (\ref{assF}) and S (\ref{assS}) hold. 
The goal of this section is to identify a suitable invariant to detect whether 
the reduction of $A$ at a good prime $\vvv$ of $L$ (see \Cref{pgood}) is $\mu$-ordinary; here we assume that the residue characteristic of $\fp$ is $p$.
When $p$ is not totally split in $F$, we further assume Assumption~A (\ref{assA}).

\begin{remark}\label{rmknp}
Below, we index over primes $\fP$ of $F$ lying above $p$.
Via the multiplication of $F$ on $A$, the identification $F\otimes_\QQ \QQ_p=\prod_{\fP|p} F_\fP$ induces a decomposition of $p$-divisible groups up to isogeny: $$A[p^\infty]\sim\bigoplus_{\fP|p} A[\fP^\infty].$$  
Let $\nu$ denote the Newton polygon $\nu(A_{\FF})$ of the reduction of $A$ modulo $\fp$ defined in \S\ref{SmuordNP}.
Thus
$\nu=\sum_{\fP|p}\nu_\fP$, where 
$\nu_\fP$ denotes the Newton polygon of the $p$-divisible group $A[\fP^\infty]$ modulo $\fp$.\footnote{Here modulo $\fp$ is a slight abuse of notation. More precisely, to have such a decomposition of $p$-divisible groups, we pass to $EL$ to have all $F$-actions defined (see \Cref{ACimpliesAD}). Therefore, we actually pick a prime $\overline{\fp}$ in $EL$ over $\fp$ and consider the reduction modulo $\overline{\fp}$; different choices of $\overline{\fp}$ permute the values of $\nu_\fP$ and their sum $\nu$ remains the same.\label{fn1}}
Furthermore, the multiplicities of the slopes of $\nu_\fP$ are divisible by the inertia degree of $\fP/p$ in $F/\QQ$; and for any prime $\fP$, if $\fP=\fP^*$, then
 the polygon $\nu_\fP=\nu_{\fP^*}$ is symmetric; and if $\fP\neq \fP^*$, then $\nu_\fP$ is dual to $\nu_{\fP^*}$  (\cite{Rapoport-Richartz}; recall that $*$ denotes complex conjugation on $F$, which induces a natural involution on primes of $F$) and the dual of a Newton polygon is found by replacing all slopes $\lambda$ by $1-\lambda$.
 By \cite[Theorem 4.2]{Rapoport-Richartz},  $\nu_\fP\leq \mu_\fP$,
 for each prime $\fP$, where $\mu_\fP$ is the lowest possible polygon for the $\fP$-piece defined by Rapoport and Richartz explicitly. The $\mu$-ordinary Newton polygon $\mu$ in \S\ref{SmuordNP} is given by $\mu=\sum_{\fP|p}\mu_\fP$.
\end{remark}

\subsection{The $\mu$-ordinary Newton polygon} 

\begin{lemma}\label{fieldK}
There exists a unique subfield $K\subset F$ such that the prime $p$ is totally split in $K/\QQ$ and every prime of $K$ above $p$ is inert in $F/K$.
\end{lemma}

\begin{proof}
When $p$ is totally split in $F$, we take $K=F$. Otherwise, we have Assumption A that $F/\QQ$ is abelian and then $K$ is the field defined in \Cref{defK} below.
\end{proof}

\begin{notation}\label{defK}
For a prime $\fP$ as above, let $f$ denote the inertia degree of $\fP$ over $\QQ$.
When $F/\QQ$ is abelian (Assumption A), let $D \subset \Gal(F/\QQ)$ be the decomposition group of $\fP$; and
let $K =F^D$ be the fixed field; in this case, $f=[F:K]$. When $p$ is totally split in $F$ (without Assumption A), we use $D$ to denote the identity group and $K=F, f=1$.
\end{notation}

The group $D$ and the field $K$ depend on $p$, although we omit $p$ from the notation. 
Recall that $F_0$ denotes the maximal totally real subfield of $F$.
Note that $K\subseteq F_0$ if and only if $\fP=\fP^*$, for each prime $\fP$ of $F$ above $p$.

Fix a suitable prime $\overline{\fP}$ above $p$ in $\bar{\QQ} \subset \CC$.\footnote{As in footnote~\ref{fn1}, different choices of $\overline{\fP}$ permute the values of $\mu_\fP$ computed in \Cref{muord} and thus their sum $\mu$ is independent of the choice. However, we would like to make a suitable choice so that $\mu_\fP$ is indeed the $\mu$-ordinary Newton polygon of the $\fP$ piece modulo $
\overline{\fp}$. To do so, we note that $\overline{\fp}$ gives an embedding of the reflex field $E\hookrightarrow \overline{\QQ}_p$; moreover, we have a natural embedding $E\hookrightarrow \Lbar=\overline{\QQ}$ given in the discussion right after Assumption D; we pick $\overline{\fP}$ such that its corresponding embedding $\overline{\QQ}\hookrightarrow \overline{\QQ}_p$ is compatible with the above two embeddings of $E$. Note that different choices of $\overline{\fP}$ satisfying the above condition provide the same formula in \Cref{muord}. \label{fn2}}
To each $\sigma:F \hookrightarrow \CC$, we associate a prime $\fP_\sigma$ of $F$ above $p$, by $\fP_\sigma:=\sigma^{-1}(\overline{\fP})$. By definition, $\fP^*_\sigma=\fP_{\sigma^*}$.

\begin{notation}
For each prime $\fP$ of $F$ above $p$, let $T_\fP:=\{\sigma\in {\rm Hom}_\QQ(F,\CC) \mid \fP_\sigma= \fP\}$. 
\end{notation}

Note that $T_\fP$ is a coset of $D$ and $T^*_{\fP}=T_{\fP^*}$.
Recall from Assumption~S (\ref{assS}) that $\Phi$ is a CM-type of $F$ and $\sigma_1\in\Phi$ satisfies $\cf(\sigma_1)=1$.

\begin{notation} \label{apnotation}
Denote $\att_\fP=|\Phi^*\cap T_\fP|$.
We write $\fP_1=\fP_{\sigma_1}$, $\att_1=\att_{\fP_1}$ and $\att'_1=\att_{\fP^*_1}$.
\end{notation}

By definition, $\att_\fP+\att_{\fP^*}=|T_\fP|=f$.
We give a combinatorial description of the 
$\mu$-ordinary Newton polygon $\mu_\fP$ for every $\fP|p$.

\begin{proposition}\label{muord}
Let $(F,\Phi)$ be a CM type as in Assumption S (\ref{assS}); we assume Assumption A (\ref{assA}) or $p$ totally split in $F$.\footnote{In \Cref{muord}, \Cref{nup}, and \Cref{nup1}, we do not need to assume that $F/\QQ$ is abelian or $p$ totally split in $F$ for the formula of $\mu_\fP$, once we replace the condition $K\subseteq F_0$ (resp. $K\not\subseteq F_0$) by $\fP=\fP^*$ (resp. $\fP\neq \fP^*$). Indeed Moonen's formula cited in the proof works in general for any CM field $F$ and any signature type.}
Let $p$ be a good prime as in Definition~\ref{pgood}.  
Notation as in~\ref{setting}, \ref{defK}, \ref{apnotation}.

\begin{enumerate}
\item Assume $K\subseteq F_0$.
Then
    \begin{enumerate}
        \item if $\fP\neq \fP_1$, then all slopes of $\mu_\fP$ equal $1/2$, and the slope $1/2$ occurs with multiplicity $nf$; and
        \item the slopes of $\mu_{\fP_1}$ equal $1/2-1/f$ (with multiplicity $f$), $1/2$ (with multiplicity $(n-2)f$), and $1/2+1/f$ (with multiplicity $f$).
    \end{enumerate}

\item Assume $K\not\subseteq F_0$. 
Then
    \begin{enumerate}
        \item 
        if $\fP\neq \fP_1,\fP_1^*$, then all slopes of $\mu_\fP$ equal $\att_\fP/f$, and the slope $\att_\fP/f$ occurs with multiplicity $nf$;
        \item 
        the slopes of $\mu_{\fP_1}$ equal $\att_1/f$ (with multiplicity $(n-1)f$), and $(\att_1+1)/f$ (with multiplicity $f$); and
        \item the slopes of $\mu_{\fP_1^*}$ equal $(\att'_1-1)/f$ (with multiplicity $f$), and $\att'_1/f$ (here with multiplicity $(n-1)f$).
    \end{enumerate}

\end{enumerate}
\end{proposition}

In particular, if $p$ is totally split in $F$, then the $\mu$-ordinary Newton polygon $\mu$ is ordinary.  
Conversely, if $[F:\QQ] >2$ or $n\geq 3$ and Assumptions F and S (\ref{assF}, \ref{assS}) hold, then $F$ embeds into the reflex field $E$ via $\sigma_1$ (and recall that $E$ is contained in the Galois closure of $F$); 
in this situation, by \cite[Theorem 1]{wedhorn} if the $\mu$-ordinary Newton polygon is ordinary, then $p$ is totally split in $F$.

\begin{proof}
Let $\fP$ be a prime of $F$ dividing $p$. Note that $f=|T_\fP|$, for all $\fP|p$.
By \cite[Section 1.2.5 and Lemma 1.3.4]{moonen04}, the slopes of the polygon $\mu_\fP$ are $0\leq a_1 \leq \cdots \leq a_n $, where 
$a_j \in \bQ$ is given by
$a_j:=\#\{\sigma\in T_\fP\mid \cf(\sigma)>n-j \}/f$,  for $1\leq j\leq n$, and each $a_j$ has multiplicity $f$. The desired slopes description for $\mu_\fP$ follows from direct computation of the formula and the following observations.

{\bf Case (1):} All primes of $F_0$ above $p$ are inert in $F/F_0$.\footnote{For the general case without the assumption that $F/\QQ$ is abelian, we only consider primes in $F_0$ below $T_\fP$ for $\mu_\fP$ and the same argument holds; same for the discussion in Case (2).} Thus
$\fP^*=\fP$
for any prime $\fP$ of $F$, and hence $\sigma\in T_\fP$ if and only if $\sigma^*\in T_\fP$.

{\bf Case (2):} All primes of $F_0$ above $p$ are split in $F/F_0$. Thus $\fP^*\neq\fP$ for any prime $\fP$ of $F$, and hence $\sigma\in T_\fP$ if and only if $\sigma^*\not\in T_\fP$. 
\end{proof}
 
\begin{lemma}\label{nup}
Notation and assumptions as in Proposition~\ref{muord} and Remark~\ref{rmknp}.
Then 
\begin{itemize}
\item $\nu_\fP=\mu_\fP$, for any $\fP\neq \fP_1,\fP_1^*$;

\item $\nu<\mu$ if and only if $\nu_{\fP^*_1}<\mu_{\fP^*_1}$.
\end{itemize} 
\end{lemma}

\begin{proof}
When $\fP\neq \fP_1,\fP_1^*$, by Proposition \ref{muord}, the polygon $\mu_\fP$ is isoclinic; we deduce $\nu_\fP=\mu_\fP$ from $\nu_\fP\leq \mu_\fP$.
Since $\nu_{\fP_1}=\nu_{\fP_1^*}$ when $\fP_1=\fP_1^*$ and $\nu_{\fP_1}$ is dual to $\nu_{\fP_1^*}$ when $\fP\neq \fP^*$, we conclude
that $\nu_{\fP_1}=\mu_{\fP_1}$ if and only if $\nu_{\fP^*_1}=\mu_{\fP^*_1}$.
Therefore, $\nu<\mu$ if and only if $\nu_{\fP^*_1}<\mu_{\fP^*_1}$.
\end{proof}

\begin{proposition}\label{nup1}
Notation and assumptions as in Notation \ref{apnotation} and Proposition~\ref{muord}.
\begin{enumerate}
\item  Assume $K\subseteq F_0$.
Then $\nu=\mu$ if and only if the smallest slope of $\nu_{\fP^*_1}$ is equal to $1/2-1/f$.
\item Assume $K\not\subseteq F_0$.
Then $\nu=\mu$ if and only if the smallest slope of $\nu_{\fP^*_1}$ is  equal  to $(\att'_1-1)/f$.
\end{enumerate}

\end{proposition}

\begin{proof}
By Lemma~\ref{nup}, $\nu=\mu$ if and only if $\nu_{\fP_1^*}= \mu_{\fP_1^*}$. 
The multiplicities of the slopes of $\nu_{\fP^*_1}$ are divisible by $f=[F:K]$, by Proposition~\ref{muord}.
This completes the proof
since the smallest slope of $\mu_{\fP_1^*}$ has multiplicity $f$, and $\mu_{\fP_1^*}$
has at most $2$ slopes of magnitude $\le \frac{1}{2}$.
\end{proof}

\subsection{Trace of Frobenius}\label{sec_Frob}

We now introduce the invariant we use to detect $\mu$-ordinary reduction. 
It suffices to restrict to a density $1$ set of primes of $L$ satisfying the following assumption.

\begin{assumption}\label{ass3}
For $\vvv$ a prime of $L$, we assume \begin{itemize} 
\item[(P1)] $\vvv$ is good, as defined in Definition \ref{pgood};
\item[(P2)] $\vvv$ has degree $1$; that is, $[L_{\vvv}:\QQ_p]=1$, for $p$ the rational prime lying below $\vvv$.
\end{itemize}
\end{assumption}

We follow notation from \S~\ref{SmuordNP}. The residue field $\FF$ of $\fp$ is now $\FF_p$ and $W(\FF)=\ZZ_p$. 
Define $V_p:=H^1_{\rm cris}(A_{\FF}/\ZZ_p)\otimes_{\ZZ_p} \bQ_p$ and recall that $\varphi$ denote the Frobenius action. By Assumption~\ref{ass3}(P2), $\varphi$ is a $\bQ_p$-linear map on $V_p$.\footnote{Of course $\varphi$ depends on $\fp$, but we drop the dependence on $\fp$ in the notation since we only use crystalline cohomology to deduce properties of the Frobenius actions on \'etale cohomology and we will only use the \'etale Frobenius in later sections with $\fp$ varying (thus we keep $\fp$ in the notation of \'etale Frobenius).}
Following notation from \S~\ref{sec_MTA}, pick $\ell \neq p$ and let $\Fr_{\vvv}$ denote a Frobenius element in $\Gal(\Lbar /L)$; thus $\rho_{A,\ell}(\Fr_\fp) \in G_{A,\ell} \subset \GL(V_\ell)$, where $V_\ell=H^1_{\text{\'et}}(A_{\Lbar},\QQ_\ell)$.\footnote{More precisely, for the entire paper, we only need to consider $\rho_{A,\ell}(\Gal(\Lbar /L))$; in other words, we consider a Galois extension $L'/L$ with $\Gal(L'/L)\cong \rho_{A,\ell}(\Gal(\Lbar /L))$, which is unramified away from finitely many primes ($\ell$ and primes of bad reduction for $A$). Thus for a good prime $\fp$, if we pick a prime of $L'$ above it, we obtain a Frobenius elememt in $\Gal(L'/L)$ and we use $\Fr_\fp$ to denote any lift of the Frobenius to $\Gal(\Lbar/L)$. Note that for given $\fp$, $\rho_{A,\ell}(\Fr_\fp)$ is well-defined up to conjugacy by elements in $G_{A,\ell}$ and this ambiguity in conjugacy does not affect the discussion below.} 

Since $F=\End^0(A_\Lbar)$, we then have $F\subset \End_{G_{A,\ell}^\circ}(V_\ell)$ (this notation means $\bQ_\ell$-linear maps on $V_\ell$ that commute with the $G_{A,\ell}^\circ$-action) and $F\subset \End(V_p \otimes_{\bQ_p} \bQ_p^{\rm un})$ compatible with the $\sigma_p$-linear map $\varphi$ (here $\bQ_p^{\rm un}$ denotes the maximal unramified extension of $\bQ_p$ and $\sigma_p$ denotes the absolute Frobenius action; by assumption on $p$, it is unramified in $EL$ and by the running assumption of this section and \Cref{ACimpliesAD}, we have all endomorphisms of $A_{L_\fp}$ are defined over $\bQ_p^{\rm un}$ and thus $F\subset \End(H^1_{\rm dR}(A_{\bQ_p^{\rm un}}/\bQ_p^{\rm un}))=\End(V_p \otimes_{\bQ_p} \bQ_p^{\rm un})$ via the canonical comparison between de Rham and crystalline cohomologies).

\begin{lemma}\label{frp1}
Recall notation from \S~\ref{sec2}: consider an abelian variety $A$ defined over a number field $L$, with $F= \End^0(A_\Lbar)$ a field; the relative dimension of $A$ is denoted by $n$.
Let $\vvv$ be a prime of $L$ satisfying Assumption~\ref{ass3}.
For $p$ the rational prime lying below $\vvv$, let $K$ be the subfield of $F$ defined in Lemma~\ref{fieldK} 
and let $f=[F:K]$.  
Suppose that either Assumptions F, S, and A (\ref{assF}, \ref{assS}, \ref{assA}) are satisfied or Assumptions F and S are satisfied and $p$ is totally split in $F$.
Then $K\subset \End(V_p)$, $\varphi$  is $K\otimes_\QQ \QQ_p$-linear on $V_p$, and $\rho_{A,\ell}(\Fr_{\vvv})$  is $K\otimes_\QQ \QQ_\ell$-linear on $V_\ell$.
\end{lemma}

\begin{proof}
If $p$ is totally split in $F/\QQ$, then $p$ is totally split in $E/\QQ$ and hence $\fp$ is totally split in $EL/L$. Therefore all endomorphisms of $A_\Lbar$ are defined on $A_{L_\fp}$ and thus the assertion follows.

Otherwise, our assumptions imply $E=F$ and we prove that all endomorphisms in $\End(A_\Lbar) \cap K$ are defined over $KL$. 
To see this, note that by our assumptions and \Cref{ACimpliesAD}, we have that all endomorphisms of $A_\Lbar$ are defined over $FL$. Following the proof idea of \Cref{ACimpliesAD}, for $\gamma \in \Gal(FL/L)\hookrightarrow \Gal(F/\bQ)$, we see that the signature $\cf'$ associated to the $\gamma$-twisted $F$-action on $A_\CC$ via $\gamma \circ s \circ \gamma^{-1}$ for $s\in F\cap \End(A_\Lbar)$ is given by $\cf'(\gamma \circ \sigma)=\cf(\sigma)$; on the other hand, let $\gamma'\in \Gal(F/\bQ)$ denote the automorphism on $F$ induced by $\gamma$ and then we have $\cf'(\sigma)=\cf(\sigma \circ \gamma')$. Therefore we have $\cf(\gamma \circ \sigma \circ \gamma')=\cf(\sigma)$ and due to Assumptions F, S and A, the same argument in \Cref{ACimpliesAD} implies that $\gamma' = \gamma^{-1}$. Therefore if $\gamma \in \Gal(FL/KL)$, then the induced $\gamma'$ fixes $K\subset F$. Then by Galois descent, we deduce that the multiplication by $K \subset F$ is defined over $KL$.

By Assumption~\ref{ass3}(P2) and the definition of $K$, 
the prime $p$ splits completely in $KL$.  Thus the endomorphisms in $K$ are defined over $\QQ_p$
and thus commute with $\varphi$ and $\rho_{A,\ell}(\Fr_\vvv)$.
\end{proof}
We thus regard $V_p$ (resp. $V_\ell$) as a free $K\otimes_\QQ\QQ_p$-module (resp. $K\otimes_\bQ \bQ_\ell$) of rank $nf$.
Define 
\begin{equation} \label{EdefWp}
W_p:= \wedge^f_{K\otimes_\QQ \QQ_p} V_p, \quad W_\ell :=\wedge^f_{K\otimes_\QQ \QQ_\ell} V_\ell,
\end{equation}
which we regard as a free $K\otimes_\QQ\QQ_p$-module (resp. $K\otimes_\bQ \bQ_\ell$-module)
of rank 
$c:=\binom{nf}{f}$.

\begin{lemma}\label{frp}
With notation and assumptions as in Lemma~\ref{frp1}, we have
\begin{enumerate}
\item ${\rm tr}_{K\otimes_\QQ\QQ_p}(\varphi \mid W_p)=\tr_{K\otimes_\QQ\QQ_\ell}(\rho_{A,\ell}(\Fr_{\vvv}) \mid W_\ell)\in\mathcal{O}_K$;
\item For any archimedean place $\infty$ of $K$, we have $|{\rm tr}_{K\otimes_\QQ\QQ_\ell}(\rho_{A,\ell}(\Fr_{\vvv}) \mid W_\ell)|_\infty \leq c p^{f/2}$. 
\end{enumerate}
\end{lemma}

\begin{proof}
\begin{enumerate}
\item Recall that $V=H^1(A_\CC, \bQ)$ is equipped with a natural $F$-action, thus we may view $V$ as a $K$-vector space. Let $G^K$ denote $\GL_K(V)\cap \GSp(V)$, the subgroup of the $\bQ$-group $\GSp(V)$ of elements that commute with the $K$-action. Here $G^K$ is a connected reductive group over $\bQ$. By \Cref{frp1}, we have $\varphi \in G^K(\bQ_p)$ and $\rho_{A,\ell}(\Fr_\fp)\in G^K(\bQ_\ell)$.\footnote{We use the $K$-action and polarization to cutout reductive groups in $\GL(V_p)$ and $\GL(V_\ell)$; both groups are isomorphic to $G^K_{\bQ_p}$ and $G^K_{\bQ_\ell}$ by \cite[Proposition 1.3.7]{Kisin17} and the canonical Betti-\'etale comparison. Admittedly, since we are working with PEL type Shimura varieties in the proof of this lemma, we may not need the full machinery in Kisin's work for our purpose; nevertheless, Kisin's work provides clean statements for what we need and hence we cite his work here as our main reference.} 
By applying \cite[Corollary~2.3.1]{Kisin17} to the PEL type Shimura variety associated to $G^K$ (and the Hodge cocharacter of $A$; thus $A$ is naturally a point on this Shimura variety), 
we have that there exists $x\in G^K(\bQ)$ such that $x$ is conjugate to $\rho_{A,\ell}(\Fr_\fp)$ (resp.\ $\varphi$) by an element in $G^K(\bQ_\ell)$ (resp.\ $G^K(\overline{\bQ}_p)$).
Define $W:=\wedge^f_K V$, which is a $K$-vector space of dimension $c$. Therefore 
${\rm tr}_{K\otimes_\QQ\QQ_p}(\varphi \mid W_p)=\tr_K(x \mid W) = \tr_{K\otimes_\bQ \bQ_\ell}(\rho_{A,\ell}(\Fr_{\vvv}) \mid W_\ell)$ and by definition $\tr_K(x \mid W)\in K$. 
By the Weil conjectures, all Frobenius eigenvalues on $V_p, V_\ell$ are algebraic integers and thus we deduce that ${\rm tr}_{K\otimes_\QQ\QQ_p}(\varphi \mid W_p)=\tr_{K\otimes_\QQ\QQ_\ell}(\rho_{A,\ell}(\Fr_{\vvv}) \mid W_\ell)\in\cO_K$. 

\item  
For any rational prime $\ell$, with $\ell\neq p$, the Galois $\QQ_\ell$-representation $W_\ell$ is a subquotient of the Galois representation $H^f_{\text{\'et}}(A,\QQ_\ell)=\wedge^f_{\QQ_\ell} H^1_{\text{\'et}}(A,\QQ_\ell)$.
Thus the result follows from the Weil conjectures.  \qedhere
\end{enumerate}
\end{proof}

With notation and assumptions as in  Lemma~\ref{frp}, we define the invariant
\begin{equation}\label{av}
a_{\vvv}=N_{K/\QQ}({\rm tr}_{K\otimes_\QQ\QQ_p}(\varphi \mid W_p))=N_{K/\bQ}(\tr_{K\otimes_\QQ\QQ_\ell}(\rho_{A,\ell}(\Fr_{\vvv}) \mid W_\ell))\in \ZZ.
\end{equation}

Recall $d=[F_0:\QQ]$, $f=[F:K]$, and set $r :=[K:\QQ]= 2d/f$. 

\begin{proposition}\label{key}
With notation and assumptions as in Lemma~\ref{frp1}:  
let $\vvv$ be a prime of $L$ that satisfies Assumption~\ref{ass3},
let  $a_{\vvv}\in\ZZ$ be as in \eqref{av}, 
then
\begin{enumerate}
\item $p^{d-1}$ divides $a_{\vvv}$; 
\item $|p^{-d+1} a_{\vvv}|_\infty\leq C p$, for $C=\binom{nf}{f}^r$;
\item if $A$ is not $\mu$-ordinary at $\vvv$, then $p^d$ divides $a_{\vvv}$.
\end{enumerate}
\end{proposition}

\begin{proof}
Let $\Gamma = \Hom_\bQ(F,\Qbar)$ and $D = \Gal(F/K)$. 
By a choice of an embedding $F\subset \CC$, 
we identify $\Gamma$ with ${\rm Hom}_\QQ(F,\CC)$ and $\Gamma/D$ with ${\rm Hom}_\QQ(K,\CC)$.
We also identified $\Gamma/D$ with the set of primes $\fP$ of $F$ above $p$; 
for each $\fP$, we choose a representative $\sigma_\fP \in \Gamma$ in the coset 
$T_\fP$ as in Notation \ref{apnotation}.   
We label the representatives of the cosets of $D$ as $\sigma_1, \ldots, \sigma_r$. (Here $\sigma_1$ is the same $\sigma_1$ in Assumption S.)

By Lemma~\ref{frp1}, $\varphi$ is a $K\otimes_\QQ \QQ_p$-linear map on the free $K\otimes_\bQ \bQ_p$-module $V_p$ of rank $nf$. 
Let $\{\lambda_i \mid i\in\{1,\dots, nf\} \}$ denote the set of eigenvalues 
of $\varphi$. 
Then the eigenvalues of $\varphi$ on $W_p$ are $$\gamma_I=\prod_{i\in I} \lambda_i, \text{ for every subset }I\subset \{1, \dots , nf\} \text{ of size } f.$$ 
Hence, 
$a_{\vvv}$ is a sum in $\mathcal{O}_\Qbar$ of elements of the form $\sigma_1(\gamma_{I_1})\cdots \sigma_r(\gamma_{I_r})$, 
for $I_j$ a subset of $ \{1, \dots , nf\}$ of size $f$, for $j=1,\dots ,r$. Note that each $\sigma_j(\gamma_I)$ is in $\Qbar \subset \CC$ and in the following, we consider its $p$-adic valuations via the embedding $\Qbar \subset \Qbar_p$ chosen in \Cref{fn2}.

Let  $I,J$ denote subsets of $ \{1, \dots , nf\}$ of size $f$. We deduce the following statements from \Cref{muord}, Lemma~\ref{nup} and Proposition~\ref{nup1}.

{\bf Case (1)}: Assume $K\subseteq F_0$.
Note that $f$ is even, and that  $\sigma_1^*$ is in the same coset as $\sigma_1$.

\begin{itemize}
  
\item For any $I$ and any $2\leq j\leq r$, the number $\sigma_j(\gamma_I)$ is divisible by $p^{f/2}$.

\item For any $I$, the number $\sigma_1(\gamma_I)$ is divisible by $p^{f/2-1}$. 

\item 
If $A$ is not $\mu$-ordinary, then $p^{1-f/2}\sigma_1(\gamma_{I})$ is not a unit in $\mathcal{O}_{\Qbar_p}$, for each $I$.

\end{itemize}

{\bf Case (2)}: Assume $K\not\subseteq F_0$.  Note that $r$ is even.
Without loss of generality, we may choose representatives satisfying $\sigma_i^*=\sigma_{r-i}$. 
Recall $\att_\fP+\att_{\fP^*}=f$, for all $\fP \mid p$; 
in particular, $\att_1+\att'_1=f$.

\begin{itemize}

\item For any $I, J$ and any $2\leq j\leq r-1$, 
the number $\sigma_j(\gamma_I)\sigma_{r-j} (\gamma_J)$ is divisible by $p^{f}$. 

\item For any $I, J$, the number $\sigma_1(\gamma_I)\sigma_r(\gamma_J)$ is divisible by $p^{f-1}$. 

\item If $A$ is not $\mu$-ordinary, then $p^{1-f}\sigma_1(\gamma_{I})\sigma_r(\gamma_J)$ is not a unit in $\mathcal{O}_{\Qbar_p}$, for each pair $I, J$.
 
\end{itemize}

We deduce that $a_{\vvv}$ is divisible by $p^{d-1}$ (recall $fr=2d$). By Lemma~\ref{frp}, $p^{-d+1}a_{\vvv}\in \ZZ$ and $| p^{-d+1}a_{\vvv} |_\infty\leq C p$, for $C={{nf}\choose{f}}^r$. 
\end{proof}

\section{Density of $\mu$-ordinary primes} \label{Sdensity}
In Theorems~\ref{Tdensityord} and \ref{densitysigma}, we prove the main results of the paper 
about the density of the set of primes 
for which the reduction of $A$ is $\mu$-ordinary.
The main technique is to prove that certain functions closely related to $a_\fp$ in \eqref{av} are non-constant on some connected components of the 
$\ell$-adic monodromy group.

Recall that $A$ is an abelian variety defined over a number field $L$.
We suppose Assumption F (\ref{assF}) holds and recall that 
$F=\End^0(A_\Lbar)$ is a CM field and $n=2 \dim(A)/[F:\QQ]$ is the relative dimension. 
We suppose Assumption S (\ref{assS}) holds; recall that $E$ denotes the reflex field of the associated Shimura variety. Recall the notation from Section \ref{sec_MTA}.
In particular, $G_{A,\ell}$ is the $\ell$-adic monodromy group of $A$.

In order to grasp the main ideas, the reader may consider the simpler but key case when assumptions A (\ref{assA}) and C (\ref{assC}) hold. In this case, $E=F$, $L^\conn = FL$ and the identity component of the group $G_{A,\ell}$ is $\Glcon = G_{\bQ_\ell}$. 

\subsection{A conjugacy-invariant function on the $\ell$-adic monodromy group} \label{secnonconstant}

\begin{lemma} \label{Ldiagram}
Under Assumptions F and S, there are natural group homomorphisms:

\[
\Gl \xtwoheadrightarrow{\pi_0} \pi_0(\Gl)  \xleftarrow[\simeq]{\quad\rho \quad} \Gal (L^{\rm conn}/L) \xtwoheadrightarrow{\kappa} \Gal (EL/L)  
\]
where by Definition~\ref{DefLconn} of $L^{\rm conn}$, the map $\rho_{A,\ell}$ from  \eqref{Erho} induces the isomorphism $\rho$. 
Moreover, if $F/\bQ$ is abelian, we have $\Gal (EL/L) \xhookrightarrow{\iota} \Gal (F/\QQ) $.
\end{lemma}

\begin{proof}
The definitions of $\pi_0$ and $\rho$ are clear.
By \Cref{ACimpliesAD}, $EL$ is the definition field of all endomorphisms of $A_\Lbar$ and thus $EL/L$ is Galois. The homomorphism $\kappa$ is then given by the natural embedding $EL \subset L^\conn$. 
Recall from \S\ref{sec_assumptions} that when $F/\bQ$ is abelian, we have $E=F$ and thus we obtain the natural inclusion $\iota$.
\end{proof}

For the rest of this section, in addition to Assumptions F and S, we either assume assumption A or we only work with primes of $L$ over $p$ that are totally split in $F/\bQ$ (and hence totally split in $E$). In particular, we either have $\Gal(EL/L)$ abelian or only consider $\sigma=\id \in \Gal(EL/L)$. We fix a suitable auxiliary prime $\ell$. 

\begin{definition}\label{def:Ssigma}
For each $\sigma \in \Gal(EL/L)$, let $\bS_\sigma$ be the set of primes $\vvv$ of $L$ 
satisfying Assumption~\ref{ass3} (meaning that they are good primes of degree $1$)
such that $\sigma=\Fr_{\vvv\vert EL}$ and $\fp \nmid \ell$.
Here, with some abuse of notation, we denote by $\Fr_{\vvv\mid EL}$ the Frobenius at $\fp$ in $ \Gal(EL/L)$.\footnote{More precisely, under our assumption, $\fp$ is unramified in $EL/L$ and thus defines a Frobenius element up to conjugacy in $\Gal(EL/L)$; due to our assumption $\Gal(EL/L)$ abelian or $\sigma = \id$, we do not have the ambiguity of conjugacy.} 
\end{definition}

Note that $\Fr_\fp$ defined in \S\ref{sec_Frob} maps to $\Fr_{\vvv\mid EL}$ via $\kappa\circ \rho^{-1}\circ \pi_0 \circ \rho_{A,\ell}$. The set $\bS_\sigma$ has density equal to $1/[EL:L]$.

For $\sigma \in \Gal(EL/L)$,
let $F^{\langle\sigma\rangle }:=\{\alpha\in F \mid \iota(\sigma)(\alpha)=\alpha\}$ when $F/\bQ$ is abelian and let $F^{\langle\sigma\rangle }:=F$ when $\sigma=\id$.
Given a prime $\vvv\in \bS_\sigma$, and $p$ the rational prime such that $\vvv|p$, 
then the subfield $K$ of $F$ defined in Lemma \ref{fieldK} satisfies $K=F^{\langle\sigma\rangle }$.

\begin{definition} \label{DGAell} For each $\gamma \in \Gal(L^\conn/L)$, we define ${G}_{A,\ell}^{(\gamma)}:=\pi_0^{-1}(\rho(\gamma))$.
For each $\sigma\in \Gal(EL/L)$, we define ${G}_{A,\ell}^{(\sigma)}:=\pi_0^{-1}(\rho(\kappa^{-1}(\sigma)))$.
\end{definition}

Note that $\Gl^{(\gamma)}$ is a connected component of $\Gl$ for $\gamma \in \Gal(L^\conn/L)$. Under our assumptions, ${G}_{A,\ell}^{(\sigma)}$ is stable under conjugacy by $G_{A,\ell}$ for $\sigma\in \Gal(EL/L)$.

Recall $V=H^1(A_\CC,\QQ)$, $V_\ell = H^1_{\text{\'et}}(A_{\bar{L}},\QQ_\ell)$ and  $\chi$ denotes the similitude character of the unitary similitude group $G$, as in Sections \ref{subsec: PEL Shimura datum}, \ref{sec_MTA}. 
For any subfield $K$ of $ F$, recall from Section \ref{sec_Frob} that $G^K=\GL_K\cap \GSp_{2g}$ is the centralizer of $K$ in $\GSp_{2g}=\GSp(V)$. The group $G^K$ is a reductive algebraic group containing $G$. Recall that $f=[F:K]$ and $d=[F_0: \QQ]$.

\begin{definition} \label{DWsigma}
For $\sigma\in \Gal(EL/L)$, set $K=F^{\langle\sigma\rangle}$ and let $\Ktil$ denote the Galois closure of $K/\bQ$.
Let $W'_\sigma$ denote the following algebraic representation of $G^K$ over $K$:\footnote{The representation by definition is defined over $\Ktil$ and one may argue by Galois descent that it is defined over $K$. In the end, we will pick the auxiliary $\ell$ such that $\ell$ splits completely in $F/\bQ$ and thus the base change of this representation to $\bQ_\ell$ is defined over $\bQ_\ell$.}
\begin{equation}\label{WK}
W'_{\sigma} \otimes_K \Ktil =
\chi^{-d}\otimes \left(\bigotimes_{\delta\in \Hom(K,\Ktil)}\left((\wedge^f_K V)\otimes_{K,{\delta}} \Ktil \right) \right).\end{equation}
\end{definition}

In the following, we pick the auxiliary prime $\ell$ such that $\ell$ splits completely in $F/\bQ$ and hence splits completely in $\Ktil/\bQ$; we pick an embedding $\Ktil\hookrightarrow \bQ_\ell$ and write $W'_{\sigma, \ell}=W'_\sigma\otimes_K \QQ_\ell$. Note that $W'_{\sigma, \ell}$ as a $G^K_{\bQ_\ell}$-representation is independent of the choice of $\Ktil\hookrightarrow \bQ_\ell$. Indeed, by the (canonical) Betti-\'etale comparison, we have 
\begin{equation*} \label{WKell}
    W'_{\sigma, \ell} \cong
\chi^{-d}\otimes \left(\bigotimes_{\delta\in \Hom(K,\bQ_\ell)}\left((\wedge^f_{K\otimes_{\bQ} \bQ_\ell} V_\ell)\otimes_{K\otimes_{\bQ} \bQ_\ell,{\delta}} \bQ_\ell \right) \right)
\end{equation*}
as $G^K_{\bQ_\ell} \cap G_{A,\ell}$-representations (recall that the natural $K\otimes_{\bQ} \bQ_\ell$-action on $V_\ell$ was given in \Cref{frp1}). Thus we may view $W'_{\sigma, \ell}$ as a representation of $\Gal(\Lbar/KL)$ over $\bQ_\ell$ and we have $\Fr_\fp$ acting on $W'_{\sigma, \ell}$ for $\fp \in \bS_\sigma$.
Here is an alternative description of $W'_{\sigma, \ell}$. Since $K\otimes_{\bQ} \bQ_\ell \cong \oplus_{\delta\in \Hom(K,\bQ_\ell)} \bQ_\ell$, we have a decomposition of $V_\ell$ as a $G^K_{\bQ_\ell} \cap G_{A,\ell}$-representation (equivalently a $\Gal(\Lbar/KL)$-representation):
\begin{equation} \label{EVelldelta}
 V_{\ell}=\bigoplus_{\delta\in \Hom(K,\bQ_\ell)} V_{\ell}^{\delta}, \text{ with } 
    V_{\ell}^{\delta}=V\otimes_{K\otimes_{\bQ} \bQ_\ell,\delta} \QQ_\ell .
    \end{equation}
We have
\begin{equation} \label{Welldelta}
W'_{\sigma, \ell} \cong \chi^{-d} \otimes \left(\bigotimes_{\delta\in \Hom(K,\bQ_\ell)} \wedge^f_{\bQ_\ell} V_\ell^\delta \right).
\end{equation}

\begin{definition}\label{Ninvariant}
For each $\sigma\in \Gal(EL/L)$, 
we denote by $b_\sigma(-)$ (resp.  by
$a_\sigma(-)$) the restriction to $\Gl^{(\sigma)}$ of the trace function 
$\tr (-\mid W'_{\sigma,\ell})$ (resp. 
$\tr (-\mid \chi^{d}\otimes W'_{\sigma, \ell})$) on $G^{F^{\langle \sigma \rangle}}_{\bQ_\ell}$.

Note that the map $b_\sigma$ is an algebraic morphism ${G}_{A,\ell}^{(\sigma)} \rightarrow \bA^1_{\bQ_\ell}$ and by a slight abuse of notation, we also use $b_\sigma$ to denote the map on $\bQ_\ell$-points $b_\sigma: {G}_{A,\ell}^{(\sigma)}(\bQ_\ell) \to \QQ_\ell$. Note that $b_\sigma$ is invariant under conjugation by $G_{A,\ell}$.
\end{definition}

We summarize some key properties from our discussion above.
\begin{lemma}\label{bsigma}
For each $\sigma\in \Gal(EL/L)$, 
we have $\Gl^{(\sigma)}\subseteq G^{F^{\langle \sigma \rangle}}_{\QQ_\ell}$. 
Furthermore, for each prime $\vvv\in\bS_\sigma$, we have  $\rho_{A,\ell}(\Fr_\vvv)\in {G}_{A,\ell}^{(\sigma)}$  and for $a_\vvv$ as in \eqref{av}
\begin{equation}\label{aK}
 \tr(\rho_{A,\ell}(\Fr_\vvv)\mid W'_{\sigma,\ell})=b_\sigma(\rho_{A,\ell}(\Fr_\vvv))=p^{-d}a_\vvv.
 \end{equation}
\end{lemma}
\begin{proof}
    The statements are all discussed above except for the last assertion. To deduce it, using the notation in \S\ref{sec_Frob}, we have $W'_{\sigma,\ell}\cong \chi^{-d} \otimes \left(\bigotimes_{\delta\in \Hom(K,\bQ_\ell)}\left(W_\ell\otimes_{K\otimes_{\bQ} \bQ_\ell,{\delta}} \bQ_\ell \right) \right)$. We then deduce the assertion by the definition of $a_\fp$ and the fact that $\chi(\rho_{A,\ell}(\Fr_\vvv))=p$, using that the trace on a tensor representation is given by the product.
\end{proof}

We give a lower bound for
the density of the set of primes of $L$ for which $A$ has $\mu$-ordinary reduction, by studying the conjugacy-invariant function $b_\sigma$ on $\Gl^{(\sigma)}$ defined in Definition~\ref{Ninvariant}. Note that when $n=1$, $A$ has CM and then our main results follow from Shimura--Taniyama formula; thus we assume $n\geq 2$ for the discussion below.

\begin{notation} 
When the relative dimension $n$ of $A$ is at least $2$, 
for $\sigma\in \Gal(EL/L)$, we set 
\[\eta_\sigma=\frac{|\{\gamma \in \Gal(L^{\conn}/L)\mid \kappa(\gamma)=\sigma,\, \, b_\sigma |_{\Gl^{(\gamma)} } \text{ is not constant}\}|}{|\pi_0(\Gl^{(\sigma)})|}.\]
In other words, $\eta_\sigma$ is the proportion of the number of connected components of $\Gl^{(\sigma)}$ on which $b_\sigma$ is not a constant function.
\end{notation}

The following proposition follows from Sawin's proof of \cite[Theorem 2.1]{Sawin16} adapted to our context; the proof uses Chebotarev density results for $\ell$-adic Lie groups by Serre discussed in \cite[\S\S 5.2, 6.2]{Serre12}.
\begin{proposition}\label{Psawin} Let $A$ be an abelian variety defined over a number field $L$.
Suppose $A$ satisfies Assumptions F and S, with relative dimension $n\geq 2$.
Then the lower density of the set of primes of $EL$ for which $A$ has ordinary reduction is at least $\eta_{\id}$. 
Moreover, under the extra Assumption A, the lower density\footnote{The reason that we only prove a result on lower density here and above is because \Cref{key} did not provide an if and only if condition for $\mu$-ordinary in all cases. For some certain cases when our $\mu$-ordinary condition in \Cref{key} is both necessary and sufficient (i.e., being $\mu$-ordinary implies $p^d\nmid a_\fp$), then Sawin's proof adapted in our setting shows that the density exists and equals the given formula. On the other hand, in our applications, these lower densities always equal $1$ and thus we conclude the density exists and equals $1$.} of the set of primes of $L$ for which $A$ has $\mu$-ordinary reduction is at least
\[ \frac{\sum_{\sigma \in \Gal(EL/L)}\eta_\sigma}{[EL:L]}. \]
\end{proposition}

\begin{proof}
We briefly sketch the proof based on the proof of \cite[Theorem 2.1]{Sawin16}. We focus on the second assertion. The proof of the first assertion is the same once we replace $L$ by $EL$ and note that for $\sigma=\id$, the $\mu$-ordinary Newton polygon is ordinary.

Let $\Gamma$ be the image of $\Gal(\overline{L}/L)$ in $\GSp_{2g}(\QQ_\ell)$. Note that $\Gamma$ is an $\ell$-adic Lie group and we equip $\Gamma$ with the Haar measure of total mass $1$.
By definition, $\Gamma$ is Zariski dense in the $\ell$-adic monodromy group $\Gl$. 

For each $\sigma \in \Gal(EL/L)$, we consider $\gamma \in \Gal(L^{\conn}/L)$ with $\gamma|_{EL}=\sigma$ and denote by $\Gamma^{(\gamma)}$ 
the image of the coset $\gamma \Gal(\overline{L}/L^{\conn})$ in $\GSp_{2g}(\QQ_\ell)$. Then $\Gamma^{(\gamma)}$ has measure $1/[L^{\conn}:L]$, and is Zariski dense in the connected component $G^{(\gamma)}_{A,\ell} \subseteq \Gl$.

By \Cref{bsigma} and Proposition \ref{key}(2), there exists $C>0$, independent of $\ell$,
such that  
\[b_\sigma (\rho_{A,\ell}(\Fr_\vvv))\in[-C,C], \forall \vvv \in \bS_\sigma. \]
For $n \in [-C,C] \cap \ZZ$, let $T_n$ be the 
subset of $\Gl^{(\gamma)}$ where $b_\sigma=n$.
Let $Z$ be the finite union of $T_n$ over all integers $n \in  [-C,C]$. Then $Z$ is a conjugacy-invariant closed subscheme of $\Gl^{(\gamma)}$.

By Chebotarev’s density theorem (see for instance \cite[\S\S 5.2, 6.2]{Serre12}), the density of the set of primes 
such that $\rho_{A,\ell}(\Fr_\vvv)$ is in $Z$ is equal to the Haar measure of $Z\cap \Gamma$. When $b_\sigma$ is not a constant function on $\Gl^{(\gamma)}$, since $\Gamma^{(\gamma)}$ is Zariski dense in $\Gl^{(\gamma)}$, the Haar measure of $Z\cap \Gamma$ is $0$.

By Proposition \ref{key}, if $\vvv \in \bS_\sigma$, then $A$ has $\mu$-ordinary reduction at $\vvv$ if $b_{\sigma}(\rho_{A,\ell}(\Fr_\vvv))$ is not integral. In other words, $\vvv \in \bS_\sigma$ is a prime of $\mu$-ordinary reduction for $A$ if $\rho_{A,\ell}(\Fr_\vvv)$
is not in $Z$. We deduce that the lower density of the set of primes $\vvv \in \bS_\sigma$ such that 
the reduction of $A$ modulo $\vvv$ is $\mu$-ordinary is at least the Haar measure of the union of $\Gl^{(\gamma)} \subset \Gl^{(\sigma)}$ over which $b_\sigma$ is not constant, which is exactly $\eta_\sigma/[EL:L]$.
We complete the proof by summing over all $\sigma \in \Gal(EL/L)$.
\end{proof}

\subsection{The Weyl group of the sympletic group} \label{Sweyl}
We briefly review some background.
Let $T$ be a maximal (split) torus of the sympletic group $\Sp_{2g,\Qbar}$ (over $\Qbar$ or over any algebraically closed field of characteristic $0$).
Without loss of generality, we may assume the symplectic pairing is
$\begin{bmatrix} 0 & I_g \\ -I_g & 0
\end{bmatrix}$
and we can suppose $T$ is the set of diagonal matrices of the form
$\mathrm{diag}(z_1, \ldots, z_g, z_1^{-1}, \ldots, z_g^{-1})$.
The action of $N_{\Sp_{2g}}(T)$ on $T$ by conjugation gives a permutation of the diagonal entries of $T$.

The Weyl group of $\Sp_{2g}$ is $N_{\Sp_{2g}}(T)/T$.
The action described above gives a homomorphism from $N_{\Sp_{2g}}(T)$ into $\mathfrak{S}_{2g}$, 
the symmetric group on $2g$ elements, via permutations of $\{\pm 1, \ldots, \pm g\}$.
The group $N_{\Sp_{2g}}(T)/T$ is isomorphic to the subgroup 
$\Sigma_g \subset \mathfrak{S}_{2g}$ consisting of permutations $w$ 
such that $w(-k)=-w(k)$ for $1\leq k \leq g$.
Furthermore, $\Sigma_g \simeq (\mathfrak{S}_2)^g \rtimes \mathfrak{S}_{g}$, with the standard semi-direct product structure. (For reference, one may consult \cite[Propositions~2.1, 4.1]{conradnotes}.)

In the next section, we reindex the variables in order to work with the structure over $F_0$.

\subsection{A group theoretic realization of the connected components of the $\ell$-adic monodromy group}\label{pi0}

Recall that $\pi_0(\Gl)$ is the component group of the $\ell$-adic monodromy group.
The goal of this section is to explicitly realize $\pi_0(\Gl)$ as a subgroup of an extension of some subquotient of the Weyl group of 
$\Sp_{2g}$ by a certain torus; see \Cref{Pdefphi}. In the simpler case when Assumptions A and C hold, we do not need the extra torus and the reader may focus on this case to see the key ideas of the proof. 

Under our assumption that $\ell$ splits completely in $F/\bQ$ and using the fact that $F_0$ is stable under the Rosati involution, we have the following lemma.

\begin{lemma}\label{lem: Symplectic forms} 
Recall that $V_\ell=H^1_{\text{\'et}}(A_{\bar{L}},\QQ_\ell)$.
There is a decomposition of vector spaces over $\bQ_\ell$ compatible with the $F_0$-action: 
    \[ V_{\ell}=\bigoplus_{\tau \in \Hom (F_0, \QQ_\ell)} V_{\ell}^{\tau}, \text{ with } 
    V_{\ell}^{\tau}=V_\ell \otimes_{F_0\otimes_{\bQ} \bQ_\ell,\tau} \QQ_\ell ,\]
   such that the symplectic form $\psi $ on $V_\ell$ obtained from the polarization of $A$ decomposes as a sum
    $\psi=\bigoplus_{\tau \in \Hom (F_0, \QQ_\ell)} \psi^\tau,$
    where $\psi^\tau$ is a symplectic form on $V_{\ell}^{\tau}$.
  \end{lemma}  

In the following, we identify  $\GSp(V_\ell,\psi)$ with a subgroup of $ \prod_{ \tau\in \Hom (F_0,\QQ_\ell)}\GSp (V_{\ell}^{\tau}, \psi^\tau)$.

\begin{notation}\label{basis}
Recall that $F_0$ is a totally real field of degree $d$ over $\QQ$, and that $n=g/d$. 
We label the elements of $\Hom (F_0, \QQ_\ell)$ by $\tau_1, \ldots, \tau_d$.

We choose a basis  $\{v_k\mid 1\leq k\leq 2g\}$ of $V_\ell$, which is compatible with the decomposition in Lemma \ref{lem: Symplectic forms} into $d$ subspaces, each of size $2n$.  We denote the dual basis under the inner pairing 
by $\{v_k^*\mid 1\leq k\leq 2g\}$.
Also, 
for each $\tau_i\in \Hom(F_0,\QQ_\ell)$ with $1 \leq i \leq d$, the set $\{v_{j+2n(i-1)}\mid 1\leq j\leq 2n \}$ is a basis of $V_{\ell}^{\tau_i}$ such that the symplectic form $\psi_i=\psi^{\tau_i}$ on $V_\ell^{\tau_i}$ is given by
$$\psi_i=\sum_{1\leq j\leq n} \left(v^*_{j+(2i-2)n}\wedge v^*_{j+(2i-1)n}\right),$$
and $\{v_{j+2n(i-1)}\mid 1\leq j\leq n \}$ (resp. $\{v_{j+2n(i-1)}\mid n+1\leq j\leq 2n \}$) is a basis of $V_\ell \otimes_{F\otimes_\bQ \bQ_\ell,\delta_{i,1}} \bQ_\ell$ (resp. $V_\ell \otimes_{F\otimes_{\bQ} \bQ_\ell,\delta_{i,2}} \bQ_\ell$), where $\delta_{i,k}: F \hookrightarrow \bQ_\ell$ are the two embeddings extending $\tau_i$. Such a basis exists because the $F$-action on $V$ induces an $F$-action on each $V^{\tau_i}_\ell$; as in \S\ref{subsec: PEL Shimura datum}, pick $\xi \in F$ totally imaginary, we then have a Hermitian form $\Psi_i$ on $V_\ell^{\tau_i}$, viewed as a free $F\otimes_{F_0, \tau_i} \bQ_\ell$-module such that $\psi_i=\tr_{F\otimes_{F_0, \tau_i} \bQ_\ell/F_0\otimes_{F_0,\tau_i}\bQ_\ell} \xi \Psi_i$; we may pick a basis of $V_\ell^\tau$ such that $\Psi_i$ is diagonal and $J_i$ denotes this matrix. Note that we have canonical decomposition $V^{\tau_i}_\ell=V_\ell\otimes_{F\otimes_{\bQ} \bQ_\ell,\delta_{i,1}} \bQ_\ell \oplus V_\ell \otimes_{F\otimes_\bQ \bQ_\ell,\delta_{i,2}} \bQ_\ell$ Let $\{v_{j+2n(i-1)}\mid 1\leq j\leq n \}$ (resp. $\{w_{j+2n(i-1)}\mid n+1\leq j\leq 2n \}$) be the projections of the vectors in the basis to  $V_\ell \otimes_{F\otimes_{\bQ} \bQ_\ell,\delta_{i,1}} \bQ_\ell$ (resp. $V_\ell \otimes_{F\otimes_{\bQ} \bQ_\ell,\delta_{i,2}} \bQ_\ell$); both sets are $\bQ_\ell$-bases in the corresponding spaces. Since $J_i$ is diagonal, a direct computation shows that we obtain the desired symplectic basis by setting $v_{j+2n(i-1)}=c_{i,j} w_{j+2n(i-1)}, n+1\leq j\leq 2n$, where $c_{i,j} \in \bQ_\ell$ are some non-zero constants depending on $\xi, J_i$.
\end{notation}

\begin{lemma} \label{Glcon}
Recall that $G$ denotes the unitary similitude group defined in \S\ref{subsec: PEL Shimura datum} and $\chi$ denotes the similitude character defined in \S\ref{sec_MTA}.
There are isomorphisms
\begin{equation} \label{eqGlcon}
 \GL_{1,\bQ_\ell}\times\prod_1^d \GL_{n,\bQ_\ell}\simeq G_{\bQ_\ell}\subset \GSp_{2g, \bQ_\ell},\quad \prod_1^d \GL_{n,\bQ_\ell}\simeq G^{\chi=1}_{\bQ_\ell}\subset \Sp_{2g,\bQ_\ell},
\end{equation}
where the first isomorphism is given by 
$(c, A_1, \dots, A_d) \mapsto {\rm diag}(A_1, c A^*_1, \dots, A_{d}, c A_{d}^*)$,
for $c\in \GL_1$ and $A_i \in \GL_{n,\bQ_\ell}$ for $1 \leq i \leq d$, 
where $A^*_i = J_i^{-1} (A^t_i)^{-1} J_i $ for $J_i$ as defined in \Cref{basis}.
For the second isomorphism, we omit $c$ (or take $c=1$).
\end{lemma}

\begin{proof}
The isomorphism in \eqref{eqGlcon}
is obtained using the $F$-action and our assumption that $\ell$ splits completely in $F/\bQ$.  Each $\GL_n$-component  
corresponds to a pair of conjugate embeddings of $F$ into $\bQ_\ell$. 
\end{proof}

Recall notation from \S\ref{sec_assumptions}. By \Cref{MderisGder} and \Cref{vasiu}, under Assumptions F and S, we  can identify $T^F_{\bQ_\ell}\Glcon$ (resp. $T^{F,\chi=1}_{\bQ_\ell}\Gl^{\circ,\chi=1}$) with the subgroup $G_{\bQ_\ell}$ (resp. $G_{\bQ_\ell}^{\chi=1}$) of $\GSp_{2g}$ described in \eqref{eqGlcon}. 

Let $T\subset G^{\chi=1}$ be the maximal torus generated by all diagonal matrices in $G^{\chi=1}$ with respect to the basis $\{v_k\mid 1\leq k\leq 2g\}$ in \Cref{basis}\footnote{By the proof of \Cref{Glcon}, the decomposition only holds over the Galois closure of $F/\bQ$; nevertheless, by Galois descent, the diagonal torus is indeed a torus defined over $\bQ$. However, this rationality property is not needed in our paper as we always work with $N_{\Sp_{2g}}(T)/T$ defined to be $N_{\Sp_{2g}(\Qbar)}(T(\Qbar))/T(\Qbar)$.} and thus $T_{\bQ_\ell} \subset T^{F,\chi=1}_{\bQ_\ell}\Gl^{\circ,\chi=1}$ is also a maximal torus. 
Note that $T$ has rank $g$ and thus it is a maximal torus of $\Sp_{2g}$. 
This gives a natural embedding of Weyl groups 
\[N_{T^{F,\chi=1}_{\bQ_\ell}\Gl^{\circ,\chi=1}}(T)/T = N_{G^{\chi=1}}(T)/T  \hookrightarrow N_{\Sp_{2g}}(T)/T.\]
In Section~\ref{Sweyl}, we identified $N_{\Sp_{2g}}(T)/T$ with a subgroup $\Sigma_g$ of $\mathfrak{S}_{2g}$. 

Consider the following subgroups of $N_{\Sp_{2g}}(T)/T$.

\begin{definition}\label{def_HW}
Define
\[H:= \left( N_{\Sp_{2g}}(G^{\chi=1})\cap N_{\Sp_{2g}}(T)\right)/T, \ {\mbox{and}} \ \Wc:=N_{G^{\chi=1}}(T)/T.\]
\end{definition}
Note that $H$ normalizes $\Wc$.

\begin{proposition} \label{Pdefphi}
There is a well-defined injective group homomorphism
\begin{equation}\label{phi}
\phi_1:\pi_0(\Gl) \hookrightarrow N_{\Sp_{2g,\bQ_\ell}}(T^{F,\chi=1}_{\bQ_\ell}\Gl^{\circ,\chi=1})(\Qbar_\ell)/\Gl^{\circ, \chi=1}(\Qbar_\ell),
\end{equation}
and an exact sequence in groups
\begin{equation} \label{Ephi2}
1 \rightarrow  (T^{F,\chi=1}/Z(M_A)^{\chi=1})(\Qbar_\ell) \rightarrow N_{\Sp_{2g,\bQ_\ell}}(T^{F,\chi=1}_{\bQ_\ell}\Gl^{\circ,\chi=1})(\Qbar_\ell)/\Gl^{\circ, \chi=1}(\Qbar_\ell) \xrightarrow{\phi_2} H/\Wc \rightarrow 1.
\end{equation}
where we recall that $M_A$ denotes the Mumford--Tate group of A. 
\end{proposition}

The first input of the proof is the following lemma.
\begin{lemma}[Serre; see \cite{Serre12}*{\S~8.3.4}]\label{lem_Serre}
 There is a natural isomorphism $\pi_0(\Gl)\simeq \pi_0(G_{A,\ell}^{\chi=1})$ induced by the embedding $\Gl^{\chi=1} \hookrightarrow \Gl$. In particular, $\Gl^{\circ,\chi=1}=(\Gl^{\chi=1})^\circ$.
\end{lemma}

\begin{proof}[Proof of \Cref{Pdefphi}]
We first construct $\phi_1$.
By \Cref{lem_Serre}, it suffices to define an injective homomorphism 
\[\pi_0(G_{A,\ell}^{\chi=1}) \hookrightarrow N_{\Sp_{2g,\bQ_\ell}}(T^{F,\chi=1}_{\bQ_\ell}\Gl^{\circ,\chi=1})(\Qbar_\ell)/\Gl^{\circ, \chi=1}(\Qbar_\ell).\] We observe that $\Gl^{\chi=1}$ normalizes $\Gl^{\circ,\chi=1}=(\Gl^{\chi=1})^\circ$ and that $\Gl^{\chi=1} \subset \Sp_{2g,\bQ_\ell}$; therefore, we have
\[\pi_0(\Gl^{\chi=1})=\Gl^{\chi=1}(\Qbar_\ell)/\Gl^{\circ, \chi=1}(\Qbar_\ell) \hookrightarrow N_{\Sp_{2g,\bQ_\ell}}(\Gl^{\circ,\chi=1})(\Qbar_\ell)/\Gl^{\circ, \chi=1}(\Qbar_\ell).\]

We claim that $N_{\Sp_{2g,\bQ_\ell}}(\Gl^{\circ,\chi=1})\subset N_{\Sp_{2g,\bQ_\ell}}(T^{F,\chi=1}_{\bQ_\ell}\Gl^{\circ,\chi=1})$ and thus we obtain the desired injective map $\phi_1$.
To prove this claim, note that by \Cref{MderisGder} and \Cref{vasiu}, we have the centralizer 
$C_{\Sp_{2g,\bQ_\ell}}(\Gl^{\circ,\chi=1})=T^{F,\chi=1}_{\bQ_\ell}$.
We also note that for any $g\in \Sp_{2g,\bQ_\ell}$, we have $g^{-1}C_{\Sp_{2g,\bQ_\ell}}(\Gl^{\circ,\chi=1})g=C_{\Sp_{2g,\bQ_\ell}}(g\Gl^{\circ,\chi=1}g^{-1})$.
Thus if $g\in N_{\Sp_{2g,\bQ_\ell}}(\Gl^{\circ,\chi=1})$, then $C_{\Sp_{2g,\bQ_\ell}}(\Gl^{\circ,\chi=1})=C_{\Sp_{2g, \bQ_\ell}}(g\Gl^{\circ,\chi=1}g^{-1})=g^{-1}C_{\Sp_{2g,\bQ_\ell}}(\Gl^{\circ,\chi=1})g$; in other words, $g T^{F,\chi=1}_{\bQ_\ell} g^{-1} = T^{F,\chi=1}_{\bQ_\ell}$. Therefore we obtain the desired claim.

Now we construct $\phi_2$.
For any $g\in N_{\Sp_{2g,\bQ_\ell}}(T^{F,\chi=1}_{\bQ_\ell}\Gl^{\circ,\chi=1})(\Qbar_\ell)$, note that $gTg^{-1}$ is still a maximal torus of $G^{\chi=1}_{\Qbar_\ell}$. All maximal tori over $\Qbar_\ell$ of the connected reductive algebraic group $G^{\chi=1}_{\Qbar_\ell}$ are conjugate to each other by an element in $G^{\chi=1}(\Qbar_\ell)$.
Thus we deduce that $gTg^{-1}=g_0Tg_0^{-1}$ for some $g_0 \in G^{\chi=1}(\Qbar_\ell)$. We observe that $g_0^{-1}g \in N_{\Sp_{2g}}(G^{\chi=1})(\Qbar_\ell)\cap N_{\Sp_{2g}}(T)(\Qbar_\ell)$. Note that the choice of $g_0$ is up to an element in $N_{G^{\chi=1}}(T)(\Qbar_\ell)$, thus we have a well-defined surjective group homomorphism
\[\Phi: N_{\Sp_{2g,\bQ_\ell}}(T^{F,\chi=1}_{\bQ_\ell}\Gl^{\circ,\chi=1})(\Qbar_\ell) \rightarrow N_{\Sp_{2g}}(G^{\chi=1})(\Qbar_\ell)\cap N_{\Sp_{2g}}(T)(\Qbar_\ell)/N_{G^{\chi=1}}(T)(\Qbar_\ell)=H/\Wc.\]
(Here the last equality is because Weyl groups are finite group schemes and thus $\Qbar$-points and $\Qbar_\ell$-points coincide.)
We observe that $\ker(\Phi) = G^{\chi=1}(\Qbar_\ell)= (T^{F,\chi=1}_{\bQ_\ell}\Gl^{\circ,\chi=1})(\Qbar_\ell)$. Therefore $\Phi$ induces a well-defined surjection $\phi_2: N_{\Sp_{2g,\bQ_\ell}}(T^{F,\chi=1}_{\bQ_\ell}\Gl^{\circ,\chi=1})(\Qbar_\ell)/\Gl^{\circ, \chi=1}(\Qbar_\ell) \xrightarrow{\phi_2} H/\Wc$ and 
\[\ker (\phi_2) = \ker (\Phi)/ \Gl^{\circ, \chi=1}(\Qbar_\ell) = T^{F,\chi=1}(\Qbar_\ell) / T^{F,\chi=1}(\Qbar_\ell) \cap \Gl^{\circ, \chi=1}(\Qbar_\ell) =  (T^{F,\chi=1}/Z(M_A)^{\chi=1})(\Qbar_\ell), \]
where the last equality follows from \Cref{vasiu}.
\end{proof}

\begin{lemma}\label{Lalpha}
Recall that we view $H$ as a subgroup of $N_{\Sp_{2g}}(T)/T\cong \fS_{2g}$, the permutation group acting on vectors in the basis in Notation \ref{basis}.
There is a group isomorphism  $$\zeta: H/\Wc\to (\mathfrak{S}_2)^d \rtimes \mathfrak{S}_d.$$
The inverse of $\zeta$ is given by 
$(\epsilon, \gamma)\mapsto [\alpha_\gamma \alpha_\epsilon]$,  where
$\alpha_\gamma\in \mathfrak{S}_{2g}$ is defined by
    \begin{equation}\label{alphagamma}
 \alpha_\gamma (j+2n(i-1))= j+ 2n (\gamma (i)-1), \text{ for all } 1\leq j\leq 2n, \, 1\leq i\leq d,
\end{equation}
and $\alpha_\epsilon=\alpha_{\epsilon,1}\cdots \alpha_{\epsilon,d}\in \prod_1^d\mathfrak{S}_{2n}\subseteq \mathfrak{S}_{2g}$ is defined by
    \begin{equation}\label{alphaepsilon} 
 \text{ for each }1\leq i\leq d:\quad \alpha_{\epsilon,i}=1 \text{ if }\epsilon_i=1 \text{ and }\alpha_{\epsilon,i}=\theta \text{ otherwise},
\end{equation}
with  $\theta =(1\, 2n)(2\, 2n-1)\dots (n\, n+1)\in \mathfrak{S}_{2n}$, the permutation group acting on $\{v_{j+2n(i-1)}\mid 1\leq j\leq 2n \}$. Here to see that $\alpha_\gamma, \alpha_\epsilon \in H$, we define $E_{\alpha_{\gamma}}\in \GL_{2g}$ to be the permutation matrix of $\alpha_\gamma$;
 and define $\delta_{\alpha_\epsilon} = \delta_{\alpha_{\epsilon,1}}\cdots \delta_{\alpha_{\epsilon,d}}$, with $\delta_{\alpha_{\epsilon,i}}$ being the identity matrix in $\GL_{2g}$ when $\alpha_{\epsilon,i}=1$ and $\delta_{\alpha_{\epsilon,i}}$ being the unique matrix in $\GL_{2g}$, which acts as the identity map on $V^{\tau_k}_\ell$ with $k\neq i$ and sends (the ordered set) $\{v_{j+2n(i-1)}\mid 1\leq j\leq 2n \}$ to $\{-v_{2ni}, -v_{2ni -1}, \dots, -v_{n+1+2n(i-1)}, v_{n+2n(i-1)}, \dots, v_1\}$ otherwise. We note that $E_{\alpha_\gamma}, \delta_{\alpha_\epsilon} \in \left( N_{\Sp_{2g}}(G^{\chi=1})\cap N_{\Sp_{2g}}(T)\right) (\bQ_\ell)$ thus $\alpha_\gamma, \alpha_\epsilon \in H$.
\end{lemma}

\begin{notation}\label{def:Balpha}
Given $\gamma \in \Gal(L^\conn/L)$, let $[\alpha]=\phi_2(\phi_1(\rho( \gamma)))$, 
with $\phi_1, \phi_2$ defined in \eqref{phi} and \eqref{Ephi2}. We define
\[B_{[\alpha]}:=E_{\alpha_\gamma} \delta_{\alpha_\epsilon} \in \left( N_{\Sp_{2g}}(G^{\chi=1})\cap N_{\Sp_{2g}}(T)\right) (\bQ_\ell),\] with notation as in \Cref{Lalpha} and $\zeta([\alpha])=(\epsilon, \gamma)$.
\end{notation}

\begin{proof}
By the choice of the basis in \Cref{basis}, we have that $E_{\alpha_\gamma}, \delta_{\alpha_\epsilon} \in \left(N_{\Sp_{2g}}(T)\right) (\bQ_\ell)$. By \Cref{Glcon}, we verify that $E_{\alpha_\gamma}, \delta_{\alpha_\epsilon} \in \left( N_{\Sp_{2g}}(G^{\chi=1})\right) (\bQ_\ell)$. Thus the inverse of $\zeta$ is well-defined as a map between sets. We now construct a group homomorphism $\zeta$ whose inverse matches with the above map between sets.

Note that $C_{\GL_{2g}}(G^{\chi=1})=\Res_{F/\bQ} \bG_m$. Thus any element in $N_{\Sp_{2g}}(G^{\chi=1})$ permutes the eigenspaces of the 
$F\otimes_{\bQ} {\bQ_\ell}$-action on $V_\ell$. Thus any element in $\left( N_{\Sp_{2g}}(G^{\chi=1})\cap N_{\Sp_{2g}}(T)\right) (\bQ_\ell)$ permutes the sets of bases of $V_\ell^{\tau_i}$, up to multiplication on each basis vector (i.e., up to multiplication by $T(\bQ_\ell)$); in other words, for ${1\leq i\leq 2d}$, define the set $\mathcal{B}_i=\{ j\mid n(i-1)< j\leq  ni\}$; then any element $\alpha \in H$ permutes the collection $\{\cB_i \mid 1\leq i\leq 2d\}$.

Moreover, for ${1\leq i\leq d}$, define the set $\mathcal{A}_i=\mathcal{B}_{2i-1}\sqcup \mathcal{B}_{2i}$. 
Since $\alpha \in H$ preserves the symplectic form, then $\alpha$ acts on the collection $\{\mathcal{A}_i\mid 1\leq i\leq d\}$.  This yields a surjective group homomorphism $\omega_2: H\to \mathfrak{S}_{d}$.

If $\alpha \in  \ker(\omega_2)\subset H$, then $\alpha$ acts on 
$\{\mathcal{B}_{2i-1}, \mathcal{B}_{2i}\}$ for each $1 \leq i \leq d$;
let $\omega_1(\alpha)_i$ be $1$ if this action is trivial and $-1$ otherwise.
Let $\omega_1(\alpha)=\prod_1^d \omega(\alpha)_i\in (\mathfrak{S}_2)^d$.
This yields a surjective group homomorphism 
$\omega_1: \ker(\omega_2) \to (\mathfrak{S}_{2})^d$.

The homomorphisms $\omega_1, \omega_2$ yield a surjective homomorphism 
$\Omega: H \to (\mathfrak{S}_2)^d \rtimes \mathfrak{S}_d$, with the 
standard semi-direct product structure. We can check directly that $(\epsilon, \gamma) \mapsto \alpha_\gamma \alpha_\epsilon$ is a section of $\Omega$ by checking its action on $\{\cB_i \mid 1\leq i\leq 2d\}$ and using the fact that the standard semi-direct product structure is exactly the semi-direct product structure induced by permutation action on  $\{\cB_i \mid 1\leq i\leq 2d\}$.

It remains to show that $\Wc= \ker(\Omega)$. 
By Lemma~\ref{Glcon}, $\Wc=\prod_{1}^d {\rm Im}( \Delta: \mathfrak{S}_{n} \to \mathfrak{S}_{n}\times \mathfrak{S}_{n})\subseteq H\cap \prod_1^{2d}\mathfrak{S}_{n}$, for $\Delta:\mathfrak{S}_{n} \to  \mathfrak{S}_{n}\times \mathfrak{S}_{n}$ the diagonal embedding. Hence, $\Wc\subseteq \ker(\Omega)$.
On the other hand, let $\alpha\in \ker(\Omega)$. There exists $w\in \Wc$ such that 
$\alpha_{\vert \mathcal{B}_{i}}=w_{\vert \mathcal{B}_{i}}$ for all odd $1\leq i\leq 2d$. 
Since both $\alpha, w$ preserve the symplectic form, we then conclude that 
$\alpha^{-1}w=1$, and so $\alpha \in \Wc$.  Thus $\ker(\Omega)=\Wc$ and the induced map
$\zeta:H/\Wc \to (\mathfrak{S}_2)^d \rtimes \mathfrak{S}_d$ is an isomorphism. 
\end{proof}

\subsection{The image of the Galois group}\label{sec_image_pi0}

In this section, given $\gamma\in \Gal (L^\conn/L)$, we restrict the options for 
$\phi_2(\phi_1(\rho( \gamma) ))$ in $H/\Wc$.

Given $t\in T^{F,\chi=1}(\Qbar_\ell)$ and $[\alpha] \in H/\Wc$, we define the subvariety $G_{[\alpha],t}$ of $N_{\Sp_{2g,\bQ_\ell}}(G^{\chi=1}_{\bQ_\ell})=N_{\Sp_{2g,\bQ_\ell}}(T^{F,\chi=1}_{\bQ_\ell}\Gl^{\circ,\chi=1})$ whose $\Qbar_\ell$-points are
\begin{equation}\label{Galpha}
    G_{[\alpha],t}(\Qbar_\ell)=\{B_{[\alpha]} t g\mid g\in \Gl^{\circ, \chi=1}(\Qbar_\ell)\},
\end{equation}
where $B_\alpha$ is defined in \Cref{def:Balpha}.

For any subfield $K$ of $F$, denote $H_K:=\left(N_{G^{K,\chi=1}}(G^{\chi=1})\cap N_{G^{K,\chi=1}}(T)\right)/T\subseteq H$.

\begin{lemma}\label{Lphi} 
With notation as in Lemma \ref{Ldiagram} and Proposition \ref{Pdefphi},
 for $\gamma\in \Gal (L^\conn/L)$, let $[\alpha]=\phi_2(\phi_1(\rho( \gamma) ))$ and $\sigma=\kappa(\gamma)$. 
 Then 
\begin{equation}\label{EBalpha}
[\alpha]\in H_{F^{\langle\sigma\rangle }}/\Wc\subseteq H/\Wc, \quad  B_{[\alpha]} \in N_{G^{F^{\langle\sigma \rangle},\chi=1}}(G^{\chi=1})\cap N_{G^{F^{\langle\sigma \rangle},\chi=1}}(T);
\end{equation}
and there exists $t\in T^{F,\chi=1}(\Qbar_\ell)$ such that
\[\Gl^{(\gamma),\chi=1}=G_{[\alpha],t}\]
\end{lemma}
\begin{proof}
    By \Cref{bsigma}, we have $\Gl^{(\gamma),\chi=1}\subset G^{F^{\langle\sigma\rangle },\chi=1}_{\bQ_\ell}$. Therefore, replacing $\Sp_{2g}$ by $G^{F^{\langle\sigma\rangle },\chi=1}$ in the proof of \Cref{Pdefphi}, we have $[\alpha]\in H_{F^{\langle\sigma\rangle }}/\Wc$. From the definition of $B_{[\alpha]}$, we observe that it preserves each eigenspace of the $F^{\langle\sigma\rangle }$-action and thus $B_{[\alpha]}\in G^{F^{\langle\sigma \rangle}}$. The other claims about $B_{[\alpha]}$ were proved in \Cref{Lalpha}. The last assertion follows from the proof of \Cref{Pdefphi}.
\end{proof}

\subsection{Main Results} \label{Smainresults}

We prove our main results, as an application of Propositions \ref{Psawin} and  \ref{Pdefphi}.

\begin{theorem} \label{Tdensityord}
Let $A$ be an abelian variety defined over a number field $L$, and   $F=\End^0(A_\Lbar)$.
Assume Assumptions F and S: $F$ is a CM field and the signature $\cf$ 
of the multiplication of $F$ on $A$ is simple.

Then, after a finite extension $L'$ of $L$, the density of the set of primes $\vvv$ of $L'$ such that the reduction of $A$ modulo $\vvv$ is ordinary is $1$. Moreover, when $F/\bQ$ is Galois, the statement holds with $L'= FL$ and if moreover the relative dimension $n\geq 2$, then the density of ordinary primes over $L$ is $1/[LF:L]$. 
\end{theorem}

\begin{proof}
Recall that the $\mu$-ordinary Newton polygon is ordinary at primes that split completely in $F/\QQ$. 
Thus the statement follows immediately from Case (1) in the proof of Theorem \ref{densitysigma}, for which Assumption A plays no role. 
When $F/\bQ$ is Galois, a density one set of primes of $FL$ lie over rational primes which split completely in $F/\bQ$ and it suffices to take $L'=FL$. When moreover $n\geq 2$, Assumptions $F$ and $S$ imply that $E=F$ and that the $\mu$-ordinary Newton polygon is not ordinary at primes above rational primes which do not split completely in $F/\bQ$ (see the paragraph after \Cref{muord}). Thus we obtain the exact density of ordinary reduction over $L$ is $1/[LF:L]$.
\end{proof}

\begin{theorem}\label{densitysigma}
Let $A$ be an abelian variety defined over a number field $L$.  
Assume Assumptions F, S and A. 
Then the density of the set of primes $\vvv$ of $L$ such that the 
reduction of $A$ modulo $\vvv$ is $\mu$-ordinary is $1$.
\end{theorem}

\begin{proof}
Recall that when $n=1$ (i.e., $A$ has CM), the assertion follows from the Shimura--Taniyama formula; thus we assume $n\geq 2$.
By Proposition \ref{Psawin}, it suffices to show 
that for each $\gamma\in \Gal(L^\conn/L)$, the function $b_\sigma$ in \Cref{Ninvariant} is non-constant on $\Gl^{(\gamma)}$ where $\sigma=\kappa(\gamma)$.

Fix $\sigma\in\Gal(EL/L)$, and write $K=F^{{\langle \sigma\rangle}}$. For each $\gamma$ with $\kappa(\gamma)=\sigma$, it suffices to prove that $b_\sigma$ is non-constant on $\Gl^{(\gamma),\chi=1}$. 
Recall from Definitions~\ref{DGAell} and \ref{Ninvariant}, that $b_\sigma=a_\sigma$ on $\Gl^{(\gamma),\chi=1}$ and thus we will study $a_\sigma$ for the rest of the proof. By Lemma~\ref{Lphi},
$\Gl^{(\gamma),\chi=1}=G_{[\alpha],t}$ for some $[\alpha]\in H_K/\Wc$ and $t\in T^{F,\chi=1}(\Qbar_\ell)$.
Let $T':=T\cap G^{\der}$; thus by \Cref{MderisGder} and \Cref{vasiu}, we have $T'_{\bQ_\ell} \subset \Gl^{\circ,\chi=1}$. 
In each instance below, we will verify that the function $a_\sigma(B_{[\alpha]} t - )$ is non-constant on $T'$ for all $[\alpha]\in H_K/\Wc$ and $t\in T^{F,\chi=1}(\Qbar_\ell)$, and hence $b_\sigma$ is non-constant on $\Gl^{(\gamma)}$ for all $\gamma$ with $\kappa(\gamma)=\sigma$, which finishes the proof. 

\medskip

{\bf Case (1)}: Assume $\sigma=\mathrm{id}$. Then $K= F$, $[\alpha]=[1]$, and $B_{[\alpha]}=\mathbb{I}_{2g}$. 
Recall $V_\ell^\delta$ from \eqref{EVelldelta}.  The function $a_{\mathrm{id}}(t-)=\prod_{\delta\in\Hom(F,\bQ_\ell)}(\tr (- \mid V_\ell^\delta))$ on $T'$ is not constant for any $t\in T^{F,\chi=1}(\Qbar_\ell)$. 

\medskip

{\bf Case (2)}: Assume $K\subseteq F_0$. Recall $r=[K:\QQ]$; set $f_0:=[F_0:K]$, and $m:=nf_0$.  Note $f=2f_0$.

Let $\GL_K(V_\ell)$ denote the invertible $\bQ_\ell$-linear maps on $V_\ell$ that commute with the $K$-action. Then 
\begin{equation}\label{decompGK}
    \GL_{K}(V_\ell)=\prod_{\Hom(K, \QQ_\ell)} \GL_{2m,\bQ_\ell}. 
\end{equation}
Note that for any $[\alpha]\in H_K/\Wc$, we have $B_{[\alpha]}\in G^{K,\chi=1}(\bQ_\ell)$. We write $B_{[\alpha]}=(B_{[\alpha],1}, \dots, B_{[\alpha],r})$ using the decomposition \eqref{decompGK}. We also have $t\in G^{K,\chi=1}(\Qbar_\ell)$ and we write $t=(t_1, \dots, t_r)$. 

Recall Notation \ref{basis}. For convenience, we order $\tau_i \in \Hom(F_0, \QQ_\ell)$ so that the following is satisfied:
for each $1\leq k\leq r$, if $(k-1)f_0< i\leq kf_0$ then $\tau_{i\vert K}=\tau_{1+(k-1)f_0 \vert K}$. 
We observe that with respect to the decomposition \eqref{decompGK}, $T'=\prod_{k=1}^r T_k$, where $T_k$ is a subtorus of the diagonal maximal torus of the symplectic group with symplectic form $\oplus_{i=(k-1)f_0+1}^{kf_0}\psi_i$. More precisely, $T_k$ is the subtorus such that each matrix $A_i$, for $1 \leq i \leq d$, for each element in \Cref{Glcon}, has determinant $1$.

Then by \Cref{DWsigma} and the discussion below it, including \eqref{Welldelta}, we have
\begin{equation} \label{Ecase2b}
    a_\sigma(B_{[\alpha]}t M)=\prod_{k=1}^r {\rm tr} (B_{[\alpha],k}t_k M_k\mid \wedge_{\QQ_\ell}^{2f_0} V_{\ell}^{\delta_k}),
    \end{equation}
where $M=(M_1,\dots, M_r)\in T'$ with $M_k \in T_k$ for $1\leq k \leq r$ and $\delta_k=\tau_i |_K \in \Hom(K,\bQ_\ell)$ for $(k-1)f_0< i\leq kf_0$. Thus it suffices to show that for each $k$, the function ${\rm tr} (B_{[\alpha],k}t_k M_k\mid \wedge_{\QQ_\ell}^{2f_0} V_{\ell}^{\delta_k})$ on $M_k\in T_k$ is non-constant. 

We make a heuristic remark here (not needed for the logic of the proof but hopefully clarifying the idea). Since the rest of the argument is purely combinatorial with no Galois theory needed, the proof for each factor is indeed the same as that of the case $r=1$ (i.e. $f_0=d$, and $K=\QQ$) because each $B_{[\alpha],k}$ is the $B$-matrix for some $[\alpha_k]$ corresponding to the symplectic group on $V_{\ell}^{\delta_k}$. 

Fix $k$ and consider $V_\ell^{\delta_k}$.
Let $$\mathcal{P}= \{1+2n(i-1) \mid (k-1)f_0 < i\leq kf_0\}\cup \{2ni \mid (k-1)f_0 < i\leq kf_0\}\subseteq \{1,\dots, 2g\},$$
which is a subset of size $2f_0$.
Note that $[\alpha]\in H_K/\Wc$ implies that $\alpha(\mathcal{P})=\mathcal{P}$. 
The set $\mathcal{P}$ corresponds to an eigenvector in $\wedge_{\QQ_\ell}^{2f_0} V_{\ell}^{\delta_k}$ of 
$B_{[\alpha],k}t_k M_k$, for $M_k\in T_k$.
The result follows because the eigenvalue for this eigenvector,
as a function on $T_k$, is non-constant and linearly independent of all other eigenvalues.

\medskip

{\bf Case (3)}: Assume $K\not\subseteq F_0$ and $K\neq F$. Then $[\alpha]\in H_K/\Wc$ implies that $\zeta([\alpha])=(1, \gamma)$, i.e., $[\alpha]=[\alpha_\gamma]$.

We follow the proof of case (2), beginning by
replacing $K$ with $K_0=K \cap F_0$. In this situation, $[F_0: K_0]=[F:K]=f$. 
Set $r_0=[K_0: \QQ]$ and $m=nf$; we have $ \GL_{K_0}(V_\ell)=\prod_{\Hom(K_0, \QQ_\ell)} \GL_{2m,\bQ_\ell}$ and write $B_{[\alpha]}=(B_{[\alpha],1}, \dots, B_{[\alpha],r_0})$ and $t=(t_1, \dots, t_{r_0})$; we order $\tau_i \in \Hom(F_0, \QQ_\ell)$ so that for each $1\leq k\leq r_0$, if $(k-1)f< i\leq kf$ then $\tau_{i\vert K_0}=\tau_{1+(k-1)f \vert K_0}$ and we have $T'=\prod_{k=1}^{r_0} T_k$. We also have the decomposition of $a_\sigma$:
\[a_\sigma(B_{[\alpha]}t M)=\prod_{k=1}^{r_0} \left( {\rm tr} (B_{[\alpha],k}t_k M_k\mid \wedge_{\QQ_\ell}^{f} V_{\ell}^{\delta_{k,1}}){\rm tr} (B_{[\alpha],k}t_k M_k\mid \wedge_{\QQ_\ell}^{f} V_{\ell}^{\delta_{k,2}})\right),\]
where $M=(M_1,\dots, M_{r_0})\in T'$ and $\delta_{k,1}, \delta_{k,2} \in \Hom(K,\bQ_\ell)$ are the two embeddings such that their restriction to $K_0$ is $\tau_i |_{K_0} $ for $(k-1)f< i\leq kf$. Note that by definition, $B_{[\alpha],k}, t_k, M_k$ commute with the $K$-action and $V_\ell^{\delta_{k,i}}, i=1,2$ are $K$-eigenspaces and thus $B_{[\alpha],k}, t_k, M_k$ act on $V_\ell^{\delta_{k,i}}$.

It suffices to show that for each $k$, the function ${\rm tr} (B_{[\alpha],k}t_k M_k\mid \wedge_{\QQ_\ell}^{f} V_{\ell}^{\delta_{k,1}}){\rm tr} (B_{[\alpha],k}t_k M_k\mid \wedge_{\QQ_\ell}^{f} V_{\ell}^{\delta_{k,2}})$ on $M_k\in T_k$ is non-constant. 
(Similarly with case (2), heuristically, the proof for each factor is the same as that of the case that $r_0=1$ (i.e., $f_0=d$ and $K$ is an imaginary quadratic field).)

Fix $k$ and consider $V_\ell^{\delta_k}$. 
We define (recall that $n\geq 2$)
$$\mathcal{P}_1=\{1+2n(i-1) \mid (k-1)f_0 < i\leq kf_0 \}, \quad \mathcal{P}_2=\{2+n+2n(i-1) \mid (k-1)f_0 < i\leq kf_0 \}.$$
Each $\cP_c$ for $c=1,2$ is of size $f_0$.
Since $[\alpha]=[\alpha_\gamma]$, the set $\mathcal{P}_1$ (resp. $\mathcal{P}_2$) corresponds to an eigenvector in $\wedge_{\QQ_\ell}^{f} V_{\ell}^{\delta_{k,1}}$ (resp. $\wedge_{\QQ_\ell}^{f} V_{\ell}^{\delta_{k,2}}$) of 
$B_{[\alpha],k}t_k M_k$, for $M_k\in T_k$, whose eigenvalue we denote by $\Lambda_1$ (resp. $\Lambda_2$). As functions on $T_k$, the product $\Lambda_1\Lambda_2$ is is non-constant and linearly independent of all other analogous pairwise products of eigenvalues of $B_{[\alpha],k}t_k M_k$.
We deduce that the function $a_\sigma(B_{[\alpha]}t M)$ is also non-constant.
\end{proof}

\begin{remark}\label{rmk_Sawin}
We compare Theorems \ref{Tdensityord} and \ref{densitysigma} with Sawin's theorem on K3 surfaces in \cite{Sawin}. For a K3 surface $X$ over a number field $L'$, let $E'$ denote the endomorphism field of the $\bQ$-Hodge structure given the transcendence part of $H^2(X_{\CC}, \bQ)$. The $\mu$-ordinary Newton polygon of the associated Shimura variety depends on the splitting behavior of a prime in the Galois closure of $E'$ and Sawin proved that the set of primes of $\mu$-ordinary reduction over $L'$ has density $1$; Sawin's proof generalizes to Kuga--Satake abelian varieties parametrized by GSpin Shimura varieties. If one considers the special case when $E'$ is CM (note that $E'$ is either totally real or CM), \Cref{GUtoK3} provides the link between Sawin's result and ours. More precisely, let $A_0$ be the CM abelian variety in \Cref{CMA0} and let $L''$ be a number field over which $A_0$ is defined. Then Sawin's theorem implies our theorems when $L'' \subset L$ by applying to the motive given by $\Hom(H^1(A), H^1(A_0))$ (indeed, we only need the corresponding \'etale realization); conversely, if $L''\subset L'$ and $E$ is CM and $E/\bQ$ is an abelian Galois extension, then \Cref{densitysigma} implies some special cases of Sawin's theorem. Note that a priori we may not have a desired CM abelian variety $A_0$ defined over a small number field $L''$; in such cases, we are not able to use \Cref{GUtoK3} to transfer results between the abelian varieties that we consider and the Kuga--Satake abelian varieties and K3 surfaces that Sawin considers.
\end{remark}

\begin{remark} \label{Rnotas}
The conclusion in Theorems \ref{Tdensityord} and \ref{densitysigma} still holds for abelian varieties which are products of abelian varieties satisfying the assumptions in the theorems.
\end{remark}

\begin{remark}\label{shortcut}
We remark on the proof of \Cref{densitysigma}. Even though $L^{\conn}$ may be strictly larger than $EL$ (in other words, there may be more than one connected component of $\Gl$ for an element $\sigma \in \Gal(EL/L)$),
the proof shows that $b_\sigma$ differs from component to component in an explicitly describable way, so that 
if $b_\sigma$ is non-constant on one component then it is non-constant for every 
component above $\sigma$. 
 In other words, primes corresponding to the same $\sigma$ share similar behavior. 

We provide a sketch of a potential alternative proof of \Cref{densitysigma} under extra assumptions. Assume that the $\mu$-ordinary condition in Proposition~\ref{key} is both sufficient and necessary (i.e., being $\mu$-ordinary implies $p^d\nmid a_\fp$) and assumptions A and D holds. 
If for each $\sigma \in \Gal(EL/L)$, there exists $\fp_1, \fp_2 \in \bS_\sigma$ and an abelian variety $B_{\sigma,1},B_{\sigma,2}$ over $L$ such that $\End^0(B_{\sigma,i, \Lbar})=F$ and $B_{\sigma,i}$ has $\mu$-ordinary reduction at $\fp_i$ for $i=1,2$. We then prove that for all $\sigma$, the function $b_\sigma$ is non-constant for $A$ and hence $A$ has a density $1$ of $\mu$-ordinary reduction. Indeed, for each $\sigma$, we can define the analogous function of $b_\sigma$ for $B_{\sigma,i}$. Assumption $D$ shows that there is at most one component of $\Gl$ over $\sigma$ and the assumption that $\End^0(B_{\sigma,i, \Lbar})=F$ ensures that $b_\sigma$ is the same function (under suitable identification of $V_\ell$ for $A,B_{\sigma,i}$). The assumption that $B_{\sigma,i}$ has $\mu$-ordinary reduction at $\fp_i$ and the last paragraph of the proof of \cite[Theorem 2.1]{Sawin16} show that $b_\sigma$ is non-constant.
\end{remark}

\begin{remark}
Theorem~\ref{Tdensityord} generalizes and strengthens \cite[Theorem 1]{Fite20} which is the case $g=3$ and $n=3$.
Other previous work having an overlap with Theorem~\ref{Tdensityord} includes \cite[Theorem~2]{Fite20}; he considered the case when the relative dimension is $n=2$. 
If we further assume $\End^0(A_{\Lbar})$ is a CM field and the hypotheses of Theorem~\ref{Tdensityord} applies, then Theorem~\ref{Tdensityord} gives a stronger result because the condition on the $p$-rank in \cite[Theorem 2]{Fite20} does not determine the Newton polygon (see Remark~\ref{rmk2} for an example).
\end{remark}

\section{Applications about reduction of curves}\label{SECcurves}

This section contains applications of Theorem \ref{densitysigma} for the abelian varieties 
in Moonen's special families. 
These abelian varieties are Jacobians of curves that are Galois covers of the projective line, 
with cyclic Galois group.  The meaning of the word special is that 
the image of the Hurwitz space under the Torelli map is open and dense in a component of the 
associated Deligne--Mostov Shimura variety.

We remark that other applications are possible when the 
Galois group of the cover of curves is abelian, but not cyclic.
A partial classification of the special families that occur for abelian covers 
can be found in a series of papers beginning with \cite{moonenoort}.

\subsection{Applications to Moonen's special families} \label{Smoonen}

We consider a curve $C$ in one of 
Moonen's 20 special families $M[\rr]$ for $1 \leq \rr \leq 20$ as in \cite[Table 1]{Moonen10}. 
Each of these families parametrizes branched cyclic covers of the projective line, with prescribed 
monodromy datum $(m, N, a)$.  
Here the degree of the cyclic cover $C \to {\mathbb P}^1$ is $m$, 
the number of branch points of the cover is $N$, and 
the inertia type of the cover is $a=(a_1, \ldots, a_N)$.
Let $\cf$ be the signature of the $\mu_m$-cover.

Such a cover has an affine equation of the form $y^m = \prod_{i=1}^{N-1}(x-b_i)^{a_i}$, 
where the branch locus is $\{b_i \mid 1 \leq i \leq N-1\} \cup \{\infty\}$ (and $b_1=0$, $b_2=1$).
Every integer $1 \leq g \leq 7$ occurs as the genus of the curves in at least one of these families.

Let $\zeta_m\in\CC$ be a primitive $m$-th root of unity.  

\begin{corollary}\label{CorMoo} 
Let $C$ be a smooth projective curve defined over a number field $L$, and
let $J_C$ be its Jacobian.
Suppose $C$ is a $\mu_m$-cover of the projective line, 
branched at $N\geq 4$ points, 
in the family $M[\rr]$ for $\rr=6,8,10,11$ or $14 \leq \rr \leq 20$
as in \cite[Table 1]{Moonen10}.\footnote{We exclude 9 of Moonen's families in Corollary~\ref{CorMoo} because the result in those cases
follows from earlier work of Serre \cite{Serre81}, Pink \cite{Pink98}, and Sawin \cite{Sawin16}. 
Specifically: when $1\leq \rr\leq 5$, then $C$ has genus $g \leq 2$;
when $\rr=7,9,12,13$, then $C$ has genus $3$ or $4$ but the absolutely simple factors of $J_C$ all have dimension $1$ or $2$.

Corollary~\ref{CorMoo} yields a new result for the Moonen families $M[\kappa](m,N,a)$ as follows:

$M[6] (3, 5, (1,1,1,1,2))$;
$M[8] (4, 5, (1,1,2,2,2))$;
$M[10] ( 3, 6, (1,1,1,1,1,1))$:
$M[11] (5, 4, (1,3,3,3))$;

$M[14] (6, 5, (2,2,2,3,3))$;
$M[15] (8, 4, (2,4,5,5))$;
$M[16]  (5, 5, (2,2,2,2,2))$;
$M[17] (7, 4, (2,4,4,4))$;

$M[18] (10,4, (3,5,6,6))$;
$M[19] (9, 4, (3,5,5,5))$; and
$M[20] (12, 4, (4,6,7,7))$.}

Assume the $\QQ$-algebra $\End^0(J_{C,\Lbar})$ is generated by the image of $\mu_m$.

Then the set of primes where the reduction of $C$ is $\mu$-ordinary has density $1$.\footnote{The $\mu$-ordinary Newton polygons of all these families can be found in \cite[Section 6]{LMPT2}.}
After passing to the finite extension $L(\zeta_m)/L$, the set of primes where the reduction of $C$ is ordinary has density $1$.
\end{corollary}

By \cite[Theorem 3.6]{Moonen10}, the condition on $\mathrm{End}^0(J_{C,\Lbar})$ is true for a generic curve in the family. 

\begin{proof}
By hypothesis, the curve $C$ is a $\mu_m$-cover of ${\mathbb P}^1$, with $m \geq 3$.
The Jacobian $J_C$ is a principally polarized abelian variety defined over $L$. 
The action of $\mu_m$ on $C$ induces multiplication on $J_{C,L(\zeta_m)}$ by the group algebra $\QQ[\mu_m]$. Moreover, the $\bQ$-algebra homomorphism $\QQ[\mu_m]\rightarrow \End^0(J_{C,L(\zeta_m)})$ is compatible with the $\Gal(L(\zeta_m)/L)$-action (with the natural action on $\QQ[\mu_m]$ via $\Gal(L(\zeta_m)/L)\subset \Gal(\QQ(\zeta_m)/\QQ)$ on the left-hand side and the $\Gal(L(\zeta_m)/L)$-action on endomorphisms on the right-hand side). 
The decomposition $\QQ[\mu_m] = \prod_{0 < d \mid m} \QQ(\zeta_d)$ is compatible with the $\Gal(\QQ(\zeta_m)/\QQ)$-action and therefore the projection to each factor induces an endomorphism of $J_C$ over $L$.
Hence, by \cite[Theorem B]{kanirosen89}, there is an $L$-isogeny 
$J_C\sim\prod_{0 < d|m} A_d$, where $A_d$ is an abelian variety over $L$ with multiplication by $\QQ(\zeta_d)$. 
Furthermore, for $A_d\neq\{0\}$, we have
$\mathrm{End}^0(A_{d,\Lbar})=\QQ(\zeta_d)$ because of the hypothesis that $\mathrm{End}^0(J_{C,\Lbar})$ is generated by the image of $\mu_m$.

If $d|m$ and $d >2$, then $\QQ(\zeta_d)$ is a CM field and $\Gal(\QQ(\zeta_d)/\QQ)$ is abelian.  
Direct computations show that the signature of $A_d$ with respect to the multiplication of $\QQ(\zeta_d)$ is simple.
For $d=1$, $\QQ(\zeta_d)=\QQ$ and the factor $A_1$ is trivial; 
for $d=2$, $\QQ(\zeta_d)=\QQ$ and the factor $A_2$ has dimension at most $1$. 

Hence, the statement follows from Theorem~\ref{densitysigma},  
(along with earlier results for 
abelian varieties of dimension $1$ or having CM).
\end{proof}

\begin{remark}\label{rmk2}
We compare Corollary \ref{CorMoo} with work of Fit\'e \cite{Fite20} and explain the overlaps.
Consider a curve $C$ in Moonen's family $M[\rr]$ (labeled as in \cite[Table 1]{Moonen10}).
\begin{enumerate}
\item 
The result \cite[Theorem~1]{Fite20} applies when
$\rr=6,8$ (when $C$ has genus $3$) and 
$\rr=14$ (when $C$ has genus $4$ but $J_C$ is not absolutely simple).  
Corollary~\ref{CorMoo} strengthens \cite[Theorem~1]{Fite20} for these cases
by finding the exact density for each Newton polygon.

\item 
When $\rr=11, 15$
and $17 \leq \rr \leq 20$, then the absolutely simple factors of $J_C$ 
either have dimension $1$ or $2$ or satisfy the hypotheses of Fit\'e's result
\cite[Theorem~2]{Fite20}. 
Corollary~\ref{CorMoo} gives a stronger result in these cases, because
the condition of having $p$-rank at least $\lceil g/2 \rceil$ occurs
for both $\mu$-ordinary and non $\mu$-ordinary Newton polygons. It also provides the exact density for each Newton polygon.

For example, when $\rr=11$, then $g=4$ and $J_C$ has multiplication by $F=\QQ(\zeta_5)$; 
if $p \equiv 1 \bmod 5$, then the 
reduction always has $p$-rank 
at least $2$, whether or not it is $\mu$-ordinary.
In the next section, we illustrate Corollary~\ref{CorMoo} by finding the densities for this example.
\end{enumerate}
\end{remark}

\begin{example}
The family $M[10]$ consists of $\mu_3$-covers $C \to {\mathbb P}^1$ branched at $N=6$ points with inertia type $(1,1,1,1,1,1)$.  Then $C$ has genus $4$
and an equation of the form $y^3=f(x)$ where $f(x)$ is a separable polynomial of degree $5$ (or $6$).
By Corollary~\ref{CorMoo},
the set of primes $\vvv$ of $L$ 
where the reduction of $C$ at $\vvv$ is $\mu$-ordinary has density $1$.
Suppose $\QQ(\zeta_3) \not \subset L$. 
If $\vvv$ lies above a prime $p \equiv 1 \bmod 3$, then the $\mu$-ordinary 
Newton polygon is ordinary; if instead $\vvv$ lies above a prime $p \equiv 2 \bmod 3$, it has $p$-rank $2$ and slopes $0$, $1/2$, and $1$; Corollary~\ref{CorMoo} shows that the set of primes of $L$ for each of the two Newton polygons mentioned above has density $1/2$.
\end{example}

\subsection{An example} 
\label{Sm11}
We conclude by proving an example of Proposition~\ref{key}, Theorem~\ref{densitysigma}, and Corollary~\ref{CorMoo} with additional detail.
We clarify that the statements in this section are already proved from the earlier sections. We essentially specialize the general proof to this example to highlight the main ideas in the proof without distraction from technical details.

Throughout this section, let $L=\QQ$ and $F=\QQ(\zeta_5)$.
Let $A$ be the Jacobian $J_C$ of a smooth projective curve $C$ in the Moonen family $M[11]$.
The curve $C$ is
given by an affine equation
\begin{equation} \label{Emoonen11}
    y^5=x(x-1)(x-t), \text{ for some } t \not = 0,1.
\end{equation}

Suppose that $C$ is defined over $L=\QQ$; for example, this is true if $t \in \QQ$.
Assume the Jacobian $J_C$ is not a CM abelian variety; by \cite[Theorem 3.6]{Moonen10},
this is true for a generic choice of $t$.

We first check that $J_C$ satisfies Assumptions F, S, and A. 

\begin{proof}[Proof of Assumptions F, S, and A] 
There is a natural action of $\mu_5$ on $C$ over $F=\QQ(\zeta_5)$, given by 
$\zeta_5\cdot (x,y)=(x,\zeta_5 y)$.
It induces an inclusion $F \hookrightarrow \End^0(J_{C,F})$. 
Because the family $M[11]$ has dimension 1, 
the hypothesis that $J_C$ is not a CM abelian variety implies that $\mathrm{End}^0(J_{C,\Qbar}) = F$.
Note that $F$ is a CM field and that $\Gal(F/\QQ)$ is abelian.

For $k=1, \dots, 4$, we denote by $\sigma_k: F \hookrightarrow \CC$ the embedding 
given by $\zeta_5\to e^{4\pi i k/5}$. By \cite[Lemma 2.7 in Section 3.2]{Moonen10},
along with \cite[Lemma 2.1 and Corollary 6.4]{LMPT2},
the signature of the multiplication by $F$ on $J_C$ is $\cf=(1,2,0,1)$. 
That is, $\cf(\sigma_1) =1$, $\cf(\sigma_2)=2$, $\cf(\sigma_3)=0$, and $\cf(\sigma_4)=1$. 
\end{proof}
 
Let $p\ne 5$ be a rational prime. 
Then $p$ is a good prime following Definition \ref{pgood}. 
The Newton polygon of the reduction of $J_C$ modulo $p$ is either the $\mu$-ordinary Newton polygon $\mu$ or the basic
Newton polygon $\nu$. Here $\mu$ and $\nu$ depend on the congruence class of $p \bmod 5$. Explicitly, 
\begin{equation}\label{munu5}
\mu=\begin{cases}  \ord^4 \text{ if } p\equiv 1\bmod 5\\ \ord^2\oplus \sss^2 \text{ if } p\equiv 4\bmod 5\\
(\frac{1}{4},\frac{3}{4}) \text{ if } p\equiv 2,3\bmod 5,
\end{cases} 
\text{ and }
\nu=\begin{cases}  \ord^2\oplus \sss^2 \text{ if } p\equiv 1\bmod 5\\ 
\sss^4 \text{ if } p\equiv 2,3,4\bmod 5.  \end{cases}
\end{equation}
Here, $\ord$ denotes the Newton polygon of slopes $\{0,1\}$, each with multiplicity $1$; $\sss$ is the Newton polygon of slope $\frac{1}{2}$, with multiplicity $2$, and $(\frac{1}{4},\frac{3}{4})$ is the Newton polygon of slopes $\{\frac{1}{4},\frac{3}{4} \}$, each with multiplicity $4$.

We fix a prime $\ell$ splitting completely in $F/\bQ$ and we consider $p\neq \ell$; we use $\Fr_p$ to denote the Frobenius action on $H^1_{\text{\'et}}(J_{C,\Qbar},\QQ_\ell)$.
We fix an embedding $ F\hookrightarrow \CC$ and let $|\cdot |_\CC$ denote the absolute value on $\CC$.
We refer the reader to \Cref{fn1,fn2} for the specific details in fixing suitable embeddings $\Qbar \hookrightarrow \CC$ and $\Qbar \hookrightarrow \Qbar_p$ so that we can also view $\sigma_i$ as embeddings of $F$ to $\Qbar$ and $\Qbar_p$ with compatible geometric meaning for signatures. 

We now illustrate the proof of Proposition~\ref{key} for the family $M[11]$.
We write $v_p$ for the $p$-adic valuation on $\overline{\QQ}_p$, where $v_p(p)=1$.

\begin{proof}[Proof of special case of Proposition~\ref{key}]
With some abuse of notation, let $J_C$ denote the reduction of $J_C$ modulo $p$. 
We write $V_p=H^1_{\rm cris}(J_C/\ZZ_p)\otimes_{\ZZ_p} \bQ_p$ and let $\varphi_p$ denote the crystalline Frobenius action on $V_p$. We write $V_\ell=H^1_{\text{\'et}}(J_{C,\Qbar},\QQ_\ell)$ and let $\Fr_p$ denote the Frobenius action on $V_\ell$. 

Depending on the congruence of $p\bmod 5$, we introduce an invariant $a_p\in \ZZ$ and prove: 
(a) $a_p/p\in \ZZ$;
(b) $|a_p/p^2|_\infty$ is bounded independently of $p$; 
and (c) if $J_C$ is basic, then $a_p/p^2\in \ZZ$.

\medskip

{\bf Case (I):} Assume $p \equiv 1 \bmod 5$.  
Then $p$ is totally split in $F/\QQ$ and $K=F$.

By assumption, $\varphi_p$ on $V_p$ (resp. $\Fr_p$ on $V_\ell$) is $F\otimes_\QQ \QQ_p$-linear (resp. $F\otimes_\QQ \QQ_\ell$-linear). Define $W_p:=V_p$ (resp. $W_\ell:=V_\ell$) as a free module of rank $2$ over the \'etale algebra $F\otimes_\QQ \QQ_p$ (resp. $F\otimes_\QQ \QQ_\ell$). 
Lemma~\ref{frp} (essentially \cite[Corollary~2.3.1]{Kisin17} and the Weil conjectures) implies that $\bbb_p:={\rm tr}_{F\otimes_\QQ \QQ_p}(\varphi_p\mid W_p)={\rm tr}_{F\otimes_\QQ \QQ_p}(\Fr_p\mid W_p)\in \mathcal{O}_F$, and $|\sigma_i(\bbb_p)|_\CC \leq 2 \sqrt{p}$, 
for all $1 \leq i \leq 4$.

Let $a_p=N_{F/\QQ} (\bbb_p)$. Then $a_p\in\ZZ$ and  $|a_p|_\infty \leq 16 p^2$, proving (b).

Let $\lambda_1, \lambda_2$ be the eigenvalues of $\varphi_p$ on the $F\otimes_\bQ \bQ_p$-module $V_p$; since $F\otimes_\bQ \bQ_p=\bQ_p^4$, a priori $\lambda_i\in \Qbar_p^4$ and indeed we have $\lambda_i \in \Fbar$ by \cite[Corollary~2.3.1]{Kisin17} (see the proof of Lemma~\ref{frp}). We pick any extension of $\sigma_i:F \hookrightarrow \Qbar_p$ to $\sigma_i:\Fbar \hookrightarrow \Qbar_p$ (the rest of the proof is independent of the choice).
By \Cref{muord} and \Cref{nup} (essentially \cite{Rapoport-Richartz}), $v_p(\sigma_2(\bbb_p))\geq \min\{v_p(\sigma_2(\lambda_1)),v_p(\sigma_2(\lambda_2))\}=1$ and $v_p(\sigma_3(\bbb_p))\geq 0$ independently of whether $J_C$ is $\mu$-ordinary or basic.  
When $J_C$ is $\mu$-ordinary, $\{v_p(\sigma_i(\lambda_1)), v_p(\sigma_i(\lambda_2)\}=\{0,1\}$ for $i=1,4$ and thus
$v_p(\sigma_1(\bbb_p))=0$ and $v_p(\sigma_4(\bbb_p))=0$; 
when $J_C$ is basic,  $v_p(\sigma_1(\bbb_p))\geq 1/2$ and $v_p(\sigma_4(\bbb_p))\geq 1/2$.

Hence, the product $a_p=N_{F/\QQ} (\bbb_p) = \prod_{i=1}^4 \sigma_i(\bbb_p)$ is divisible by $p$, and if $J_C$ is basic, then $a_p$ is divisible by $p^2$, proving (a) and (c).

\medskip

{\bf Case (II):} Suppose  $p \equiv 4 \bmod 5$.
Then $K=F_0 = \QQ(\sqrt{5})$ and $p$ is split in $F_0/\QQ$ and each of the two primes of $F_0$ above $p$ are inert in $F/F_0$. Write $\Gal(F_0/\QQ)=\{\tau_1, \tau_2\}$; similar to $\sigma_i$, we can also view $\tau_i$ as embeddings of $F$ to $\Qbar, \CC, \Qbar_p$.

By assumption, $\varphi_p$ on $V_p$ is not $F\otimes_\QQ \QQ_p$-linear, but is $F_0\otimes_\QQ \QQ_p$-linear; similarly $\Fr_p$ is also only $F_0\otimes_{\bQ} \bQ_\ell$-linear.
View $V_p$ and $V_\ell$ as free modules of rank $4$ over the \'etale algebras $F_0\otimes_\QQ \QQ_{p}$ and $F_0\otimes_\QQ \QQ_{\ell}$ respectively; consider $W_p :=\wedge^2_{F_0\otimes_\QQ \QQ_{p}} V_p$ and $W_\ell :=\wedge^2_{F_0\otimes_\QQ \QQ_{\ell}} V_\ell$.
By Lemma~\ref{frp}, we deduce that $\bbb_p:={\rm tr}_{F_0\otimes_\QQ \QQ_p}(\varphi_p\mid W_p)={\rm tr}_{F_0\otimes_\QQ \QQ_\ell}(\Fr_p\mid W_\ell)\in \mathcal{O}_{F_0}$, and $|\tau_i(\bbb_p)|_\CC \leq 6 {p}$, for all $1 \leq i \leq 2$. 

Let $a_p=N_{F_0/\QQ}(\bbb_p)$. Then $a_p\in\ZZ$ and $|a_p|_\infty\leq 36 p^2$, proving (b). 

Let $\lambda_1, \lambda_2, \lambda_3, \lambda_4$ be the eigenvalues of $\varphi_p$ on the $F_0\otimes_\bQ \bQ_p$-module $V_p$; similar to Case (I), we have $\lambda_i \in \Fbar$ and we pick any extension of $\tau_i:F_0 \hookrightarrow \Qbar_p$ to $\tau_i:\Fbar \hookrightarrow \Qbar_p$.
Then $\bbb_p=\sum_{\{i,j\}\subset \{1, 2,3,4\}} \lambda_{i}\lambda_{j}$, and
 \[
a_p= N_{F_0/\QQ}(\bbb_p)=\left(\sum_{\{i,j\}\subset \{1, 2,3,4\}} \tau_1(\lambda_{i})\tau_1(\lambda_{j}) \right)
\left(\sum_{\{i,j\}\subset \{1, 2,3,4\}} \tau_2(\lambda_{i})\tau_2(\lambda_{j}) \right).
\]

By \Cref{muord} and Lemma~\ref{nup},  $v_p(\tau_2(\lambda_i))=1/2$ for $1\leq i\leq4$, whether $J_C$ is $\mu$-ordinary or basic; and we may index $\lambda_i$ in a suitable way such that
when $J_C$ is $\mu$-ordinary, $v_p(\tau_1(\lambda_i))=0$ for $i=1,2$, and $v_p(\tau_1(\lambda_i))=1$ for $i=3,4$; when $J_C$ is basic,  $v_p(\tau_1(\lambda_i))=1/2$ for $1\leq i\leq 4$.
In particular, $\tau_2(\lambda_{i})\tau_2(\lambda_{j})$  is divisible by $p$ for any $\{i,j\}\subset \{1,2,3,4\}$, proving (a).

If $J_C$ is basic, then $v_p(\tau_1(\lambda_{i})\tau_1(\lambda_{j}))>0 $
for all $\{i,j\}\subset \{1,2,3,4\}$, and in turn $v_p(a_p)>1$, proving $(c)$.

\medskip

{\bf Case (III):} Suppose $p \equiv 2,3 \bmod 5$. 
Then $p$ is totally inert in $F/\QQ$. 

By assumption, $\varphi_p$ on $V_p$ is not $F\otimes_\QQ \QQ_p$-linear, nor $F_0\otimes_\QQ \QQ_p$-linear, but it is $\QQ_p$-linear. 
View $V_p$ (resp. $V_\ell$) as a vector space of dimension 8
over $\QQ_{p}$ (resp. $\bQ_\ell$). Consider $W_p=\wedge^4_{\QQ_p} V_p$ and $W_\ell=\wedge^4_{\QQ_\ell} V_\ell$. 
By Lemma~\ref{frp},  $a_p:={\rm tr}_{\QQ_p}(\varphi_p\mid W_p)= {\rm tr}_{\QQ_\ell}(\Fr_p\mid W_\ell)\in \ZZ$, and $|a_p|_\infty \leq 70 {p}^2$, proving (b).

Let $\lambda_1, \dots, \lambda_8$ be the eigenvalues of  {$\varphi_p$ on $V_p$}; they all lie in $\Fbar$. Then $\displaystyle a_p=\sum_{\{i_1,\dots ,i_4\}\subset \{1, \dots ,8\} } \prod_1^4\lambda_{i_j}$.

We pick an embedding $\Fbar \hookrightarrow \Qbar_p$.
By \Cref{muord} and Lemma~\ref{nup}, we may index $\lambda_i$ in a suitable way such that when $J_C$ is $\mu$-ordinary, $v_p(\lambda_i)=1/4$ for $1\leq i\leq 4$, and $v_p(\lambda_i)=3/4$ for $5\leq i\leq 8$; when $J_C$ is basic, $v_p(\lambda_i)=1/2$ for all $1\leq i\leq 8$.
In particular, for any
$\{i_1, i_2, i_3,i_4\} \subset \{1, \dots ,8\}$, the  monomial  $\lambda_{i_1} \lambda_{i_2} \lambda_{i_3}\lambda_{i_4}$ is divisible by $p$, proving (a). 

If $J_C$ is basic, then each monomial  $\lambda_{i_1} \lambda_{i_2} \lambda_{i_3}\lambda_{i_4}$ is divisible by $p^2$, and in turn $a_p/p$ is divisible by $p$,  
proving (c).
\end{proof}

We now illustrate the proof of Theorem \ref{densitysigma} for the special family $M[11]$. By \Cref{AssCeqv}, a direct computation shows that $M[11]$ satisfies assumption C and hence $\pi_0(\Gl)\cong \Gal(F/\bQ)$ (and we will freely use this identification in the proof).
Note that Case (3) in the proof of Theorem \ref{densitysigma} does not occur for $M[11]$, since all proper subfields of $F$ are contained in $F_0$. 

\begin{proof}[Proof of special case of Theorem \ref{densitysigma}]
For each $\sigma\in \Gal(F/\QQ)$, $ \Gl^{(\sigma)}$ denotes the corresponding connected component of $\Gl$.
We introduce a trace function $a_\sigma(-)$ on $ \Gl^{(\sigma)}$
satisfying $a_\sigma (\Fr_p)=a_p$ for all primes $p$ satisfying $\Fr_p\in \Gl^{(\sigma)}$.  Then we prove $a_\sigma(-)$ is non-constant on $(\Gl^{(\sigma),\chi=1})$. 
We conclude that 
the set of primes of $\mu$-ordinary reduction for $A$ has density $1$ 
by \Cref{Psawin}. 

To prove non-constancy, let $F^{\langle \sigma\rangle}=\{a\in F \mid \sigma(a)=a\}$. Using the observation that $\Gl^{(\sigma)} \subset G^{F^{\langle \sigma\rangle}}_{\bQ_\ell}$ normalizes $\Glcon$, we proved that $\Gl^{(\sigma),\chi=1}=G_{[\alpha]}$, for some $\alpha\in H_{F^{\langle \sigma\rangle}}/\Wc$ (see \Cref{Pdefphi}, \Cref{Lalpha}, and \Cref{Lphi}), where $G_{[\alpha]}:=G_{[\alpha],1}$ is defined in \S\ref{sec_image_pi0} (i.e., take $t=1$ since assumption C holds), and $H_{F^{\langle \sigma\rangle}}, \Wc$ are defined in \Cref{def_HW} and \S\ref{sec_image_pi0}. Our $a_\sigma$ is also defined on every $G_{[\alpha]}$ and we will verify that it is non-constant on every $G_{[\alpha]}$, which finishes the proof.

Note that $d=n=2$. We pick a basis of $V_\ell$ as in \Cref{basis} and as in \Cref{Glcon}, we have (here all algebraic groups are over $\bQ_\ell$)
\begin{equation}\label{eqn_decomp}
    \Gl^{\circ,\chi=1}\simeq \GL_2\times\GL_2 \subseteq G^{F_0,\chi=1}_{\bQ_\ell}\cong \Sp_4\times \Sp_4\subseteq \Sp_{8},
\end{equation}
given by $(A_1, A_2)\mapsto {\rm diag}(A_1, A^*_1, A_2, A_2^*)$.
Let $T$ denote the diagonal maximal torus in $\Gl^{\circ,\chi=1}$. We use $(a_1,b_1,a_2,b_2)$ to denote the element in $T$ given by $A_i=\diag(a_i,b_i)$ for $i=1,2$ and we have $A_i^*=\diag(a_i^{-1},b_i^{-1})$.

For each $\sigma$, we use $W_{\sigma, \ell}$ to denote the $W_\ell$ in each case in the proof of special case of Proposition~\ref{key} above.
\medskip

{\bf Case (1)}: Let $\sigma=\mathrm{id}$. Then  $\Fr_p\in \Glcon$ if and only if $p\equiv 1\bmod 5$, and $F^{\langle\sigma\rangle}=F$.
Define $a_{\mathrm{id}}={\rm tr}\left( - \mid \bigotimes_{{\zeta}\in \Hom(F,\QQ_\ell)} (W_{\mathrm{id},\ell}\otimes_{F\otimes_\bQ \bQ_\ell,{\zeta}} \bQ_\ell)\right)$. (Note that $F\otimes_\bQ \bQ_\ell \cong \bQ_\ell^4$ and each $\zeta$ corresponds to the projection to each coordinate.)

By definition, $H_F/\Wc=\{1\}$ thus $G_{[\alpha]}=\Gl^{\circ,\chi=1}$.
For $M=(a_1,b_1,a_2, b_2)\in T$, 
we compute
\[a_{\mathrm{id}}(M)=(a_1+b_1)(a_2+b_2)(a_1^{-1}+b_1^{-1})+(a_2^{-1}+b_2^{-1}).\] 
Then the function $a_{\mathrm{id}}(\cdot)$ is non-constant on $T$, and hence also on $\Glcon$. 

\medskip 

{\bf Case (2) part 1}: 
Let $\sigma=\star$ be the restriction of complex conjugation to $F$.  
Then  $\Fr_p\in \Gl^{(\star)}$ if and only if $p\equiv 4\bmod 5$, and $F^{\langle\sigma\rangle}=F_0$. 
Define $a_{\star}={\rm tr}\left( - \mid \bigotimes_{\tau\in \Gal(F_0/\QQ)} (W_{\star,\ell}\otimes_{F_0\otimes_\bQ \bQ_\ell,\tau} \bQ_\ell)\right)$.

We compute the subgroup $H_{F_0}/\Wc\subset H/\Wc$.
By definition, for $\alpha\in H_{F_0}$, its projection in each $\Sp_4$-factor in \eqref{eqn_decomp} normalizes the corresponding $\GL_2$-factor of $\Glcon$.
Thus we write $\alpha=(\alpha_1,\alpha_2)\in \Sigma_2 \times \Sigma_2$ (recall the Weyl group of $\Sp_{2g}$ given in \S\ref{Sweyl}).

Recall our choice of basis of $V_\ell$: $\{v_1,\dots,v_8\}$ with the symplectic forms given by $v^*_1\wedge v^*_3 + v^*_2 \wedge v^*_4, v^*_5\wedge v^*_7 + v^*_6 \wedge v^*_8$; then for the first $\Sp_4$ (similar computation holds for the second copy as well)
$$\Sigma_2=\{(1234), (4321), (12)(34), (13)(24), (14)(23), (13), (24), \mathrm{id}\}.$$
Among these, the permutations that normalize $\GL_2$ are 
\[\{(12)(34), (13)(24), (14)(23), \mathrm{id}\}\] 
and the Weyl group of $\GL_2$ is $\{(12)(34), \mathrm{id}\}$. 
We deduce that the following is a set of representatives of $H_{F_0}/\Wc$:
$$H_{F_0}/\Wc=\{\left(\mathrm{id},\mathrm{id}\right), \left((14)(23), \mathrm{id}\right), \left(\mathrm{id},(58)(67)\right), \left((14)(23), (58)(67)\right)\}.$$  

For each $[\alpha]=([\alpha_1],[\alpha_2])\in H_{F_0}/\Wc$, we define $B_{[\alpha]}=(B_{[\alpha_1]},B_{[\alpha_2]})\in \Sp_4\times \Sp_4=G^{F_0}_{\bQ_\ell}$ as follows: if $[\alpha_i]=\id$, then $B_{[\alpha_i]}:=\mathbb{I}_4$; if $[\alpha_i]\neq \id$, then $B_{[\alpha_i]}:=\begin{bmatrix}
    0 & 0 & 0 & 1\\
    0 & 0 & 1 & 0\\
    0 &-1 & 0 & 0\\
    -1& 0 & 0 & 0
\end{bmatrix}$. Then $ \Gl^{(\star)}=B_{[\alpha]}\Glcon$, which is exactly $G_{[\alpha]}$ in \Cref{Lphi}. 

Let $M=(a_1,b_1,a_2, b_2)\in T$. Then $M=(M_1,M_2)\in \Sp_4\times \Sp_4$, where 
$$M_1={\rm diag}(a_1,b_1, a_1^{-1}, b^{-1}_1)\text{ and } M_2={\rm diag}(a_2,b_2, a_2^{-1}, b^{-1}_2).$$ 
Then, for each $[\alpha]=([\alpha_1],[\alpha_2])\in H_{F_0}/\Wc$, we compute

$$a_\star(B_{[\alpha]} M)= {\rm tr}(B_{[\alpha_1]}  M_1 \mid W_{\star,\ell}^{\tau_1})\cdot {\rm tr}(B_{[\alpha_2]} M_2 \mid W_{\star,\ell}^{\tau_2}), $$
where for $[\alpha_i]=\mathrm{id}$, $${\rm tr}(B_{[\alpha_i]} M_i \mid W_{\star,\ell}^{\tau_i})= a_ib_i+a_ib_i^{-1}+a_i^{-1}b_i+a_i^{-1}b_i^{-1} +2, $$
and for $\beta\neq \id$, $${\rm tr}(B_{[\alpha_i] M_i} \mid W_{\star,\ell}^{\tau_i}))= a_ib_i^{-1}+a_i^{-1}b_i. $$

In each instance, $a_\star(B_{[\alpha]} \cdot )$ is a non-constant function on $T$, and hence $a_\sigma(\cdot)$ is non-constant on $G_{[\alpha]}$ for every $[\alpha]\in H_{F_0}/\Wc$. 

\medskip

{\bf Case (2) part 2}: Let $\sigma$ be an element of order $4$. 
Then $\rho_{A,\ell}(\Fr_p)\in\Gl^{(\sigma)}$ 
if and only if  $p\equiv 2,3 \bmod 5$, and $F^{\langle\sigma\rangle}=\QQ$.
Define $a_{\sigma}={\rm tr}\left( - \mid W_{\sigma,\ell}\right)$.

By definition, $H_\QQ=H$, and by Lemma \ref{alphagamma}, $H/\Wc$ has size $8$. We carry out the computation in one example. Consider  $[\alpha] =[(1548)(3726)] \in H/\Wc$.

Let $M=(a_1,b_1,a_2, b_2)\in T$. 
We compute
$$a_\sigma(B_{[\alpha]} M)={\rm tr}(B_{[\alpha]} M \mid W_{\sigma,\ell})= a_1b^{-1}_1a_2b^{-1}_2+a_1^{-1}b_1 a_2^{-1}b_2.$$
Then $a_\sigma (B_\alpha - ) $ is not constant on $T$, and hence the function $a_\sigma(-) $ is non-constant on $G_{[\alpha]}$.
\end{proof}


\end{document}